\RequirePackage{fix-cm}
\documentclass[twocolumn]{svjour3}          
\smartqed  

\usepackage{graphicx}
\usepackage{amsmath}
\usepackage{esint}
\usepackage{dirtytalk}
\usepackage[normalem]{ulem}
\usepackage{natbib}
\usepackage[font=scriptsize]{subfig}
\usepackage{epstopdf}
\usepackage{multirow}
\usepackage[dvipsnames]{xcolor}
\usepackage[colorlinks=true,citecolor=blue,linkcolor=ForestGreen]{hyperref}
\usepackage{mwe}
\usepackage{makecell}
\usepackage{hhline}
\usepackage{float}
\usepackage{bookmark}
\usepackage{tabu}
\usepackage{ulem} 
\usepackage{amsfonts} 
\usepackage{physics} 
\usepackage{relsize} 
\usepackage[bottom]{footmisc} 
\usepackage{enumitem} 
\setlist[itemize]{parsep=0pt, topsep=0pt, itemsep=0pt} 
\usepackage{booktabs} 
\usepackage{threeparttable} 
\usepackage{multirow} 
\usepackage{pdflscape} 
\usepackage[counterclockwise]{rotating} 
\usepackage[font=small]{caption} 
\usepackage{cuted} 
\usepackage{upgreek} 

\makeatletter 
\renewcommand*\env@matrix[1][\arraystretch]{%
	\edef\arraystretch{#1}%
	\hskip -\arraycolsep
	\let\@ifnextchar\new@ifnextchar
	\array{*\c@MaxMatrixCols c}}
\makeatother

\begin{document}
	
\title{On Mesh Deformation Techniques for Topology Optimization of Fluid-Structure Interaction Problems \thanks{Part of this work was presented at WCSMO 14.}}


\author{Mohamed Abdelhamid  \and Aleksander Czekanski }


\institute{M. Abdelhamid \at
	Department of Mechanical Engineering, Lassonde School of Engineering, York University, Toronto, Ontario, Canada \\
	\and
	A. Czekanski \at
	Department of Mechanical Engineering, Lassonde School of Engineering, York University, Toronto, Ontario, Canada \\
	\email{alex.czekanski@lassonde.yorku.ca}
}

\date{Received: date / Accepted: date}

\maketitle

\interfootnotelinepenalty=10000

\begin{abstract}
Fluid-structure interactions are a widespread phenomenon in nature. Although their numerical modeling have come a long way, the application of numerical design tools to these multiphysics problems is still lagging behind. Gradient-based optimization is the most popular approach in topology optimization currently. Hence, it's a necessity to utilize mesh deformation techniques that have continuous, smooth derivatives. In this work, we \textcolor{black}{address} mesh deformation techniques for structured, quadrilateral meshes. We discuss and comment on two legacy mesh deformation techniques; namely the spring analogy model and the linear elasticity model. In addition, we propose a new technique based on the Yeoh hyperelasticity model. We focus on mesh quality as a gateway to mesh admissibility. We propose layered selective stiffening such that the elements adjacent to the fluid-structure interface - where the bulk of the mesh distortion occurs - are stiffened in consecutive layers. The legacy and the new models are able to sustain large deformations without deprecating the mesh quality, and the results are enhanced with using layered selective stiffening.
\end{abstract}

\section{Introduction}
\label{sec:introduction}

\textbf{F}luid-\textbf{S}tructure \textbf{I}nteractions (FSI) are a common occurrence in nature. In FSI, solid structures interact with internal or surrounding fluid flows through deformation or rigid body motion. It is observed in numerous fields at different length scales; aeroelasticity where air flows over a flexible airplane wing and biomechanics where blood flows inside elastic valves/vessels \citep{Richter2017}. Although the numerical modeling of fundamental FSI problems has advanced greatly in the recent years, the application of numerical design tools such as \textbf{t}opology \textbf{o}ptimization (TO) to these problems is still in its early days. So far, the number of works that address high-fidelity \textbf{t}opology \textbf{o}ptimization of \textbf{f}luid-\textbf{s}tructure \textbf{i}nteractions (TOFSI) is still very limited compared to those that address either TO or FSI separately \citep{Lundgaard2018, Alexandersen2020}. Probably \textcolor{black}{some} of the reasons behind this lag \textcolor{black}{are} the multidisciplinary nature of TOFSI problems, which requires expertise in a wide range of technical subjects. There is also the strong nonlinear nature of FSI problems that is even exacerbated by TO resulting in numerous convergence and stability issues. An additional reason that is highly relevant to this work is the lack of proper computational tools that could be readily adopted to perform TOFSI. In TO of single physics systems, researchers are typically capable of implementing a black-box solver to handle the governing equations. This is not possible with FSI, since almost all commercial and open-source software that have FSI analysis capability treat the fluid and solid computational domains as separate entities in what is known as the \textit{separated/segregated} domain formulation. Hence, finite elements (or any other discretization unit) cannot exist simultaneously in both the solid and fluid domains, which undermines the concept of interpolation in density-based and level set TO methods. To implement a \textit{unified} domain formulation, where solid and fluid domains overlap, proprietary code has to be written to solve the FSI governing equations\footnote{While it is possible to link individual black-box solver modules to handle the FSI governing equations, it often comes at a cost in computational efficiency and functional capability. Another approach is to use binary structures to avoid interpolation, see \citep{Picelli2020}.} \citep{Yoon2010}. An essential part of this numerical platform is the mesh moving technique.

A mesh moving technique must address two critical properties: \textbf{(i)} Mesh Admissibility; overlapping and inverted elements are not permitted, and \textbf{(ii)} Mesh Quality; mainly skewness (non-orthogonality) and volume change of elements (local mesh size vanishing) \citep{Palmerio1994}. It's true that mesh admissibility is a more critical measure of failure for fluid flow solution. In addition, although mesh quality tends to have an effect in terms of degradation of solution accuracy and potential loss of positiveness (i.e., convergence issues), this effect is usually minimal compared to mesh admissibility. However, mesh quality can be considered as a gateway to mesh admissibility, meaning that maintaining high mesh quality does inherently maintain mesh admissibility. In addition, research has matured enough to fully guarantee mesh admissibility in most cases, so in this work we will focus on mesh quality as our primary goal. 

\textbf{In this work}, we discuss and propose some improvements to two legacy mesh deformation techniques for TOFSI; namely the discrete \textcolor{black}{spring analogy} model and the continuous linear elasticity model. In addition, we propose a new hyperelastic model (Yeoh) as a mesh deformation technique. We study some tuning parameters and their effect on mesh quality and computational costs. We support our findings with numerical examples where applicable. Finally, we discuss the sensitivity analysis of the new hyperelastic model in the context of a three-field FSI formulation. A few worthy remarks are in order before proceeding: \textbf{(i)} we focus mainly on mesh deformation of structured, quadrilateral meshing as it's more common in TO applications in general, \textbf{(ii)} we elect to direct our attention to mesh quality metrics as a gateway to mesh admissibility, and \textbf{(iii)} the mesh quality metrics we elect to use are: skewness to represent shape changes and change in element area to represent volume changes.

\section{Test Problems and Mesh Quality Metrics}
\label{sec:numericexmp}
\textcolor{black}{A common approach to handling mesh deformations is to model the mesh as a fictitious elastic or hyperelastic, continuous or discrete media. This way the mesh is represented by a stiffness matrix to describe its behavior, which can be linear or nonlinear depending on the material properties and/or geometric configurations. The boundary conditions of the mesh moving problem are as follows; the structural deformations at the fluid-structure interface are applied as prescribed displacements and the remaining domain boundaries are set to zero displacement.	The resulting system of equations is then solved (usually iteratively) for the nodal mesh displacements as the unknown degrees of freedom.}

\begin{figure}
	\centering
	\includegraphics[width=\linewidth]{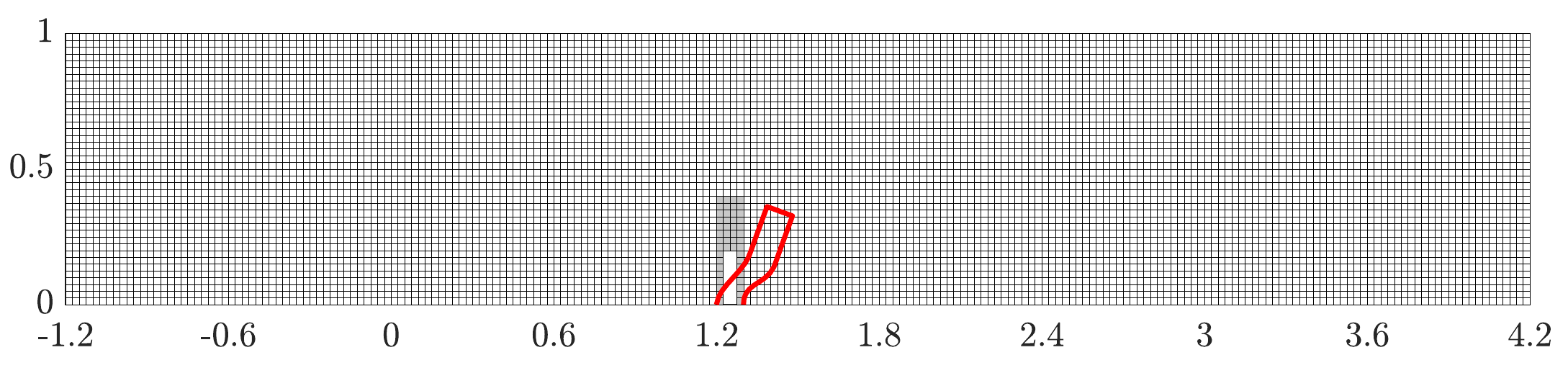}
	\caption{The original structured, quad mesh of the beam in a channel test problem. The original structure is shown in gray, while the deformed fluid-structure interface is shown in red. Dimensions are in meters unless otherwise noted.}
	\label{fig:beaminachanneltestprob}
\end{figure}

\begin{figure}
	\centering
	\includegraphics[width=\linewidth]{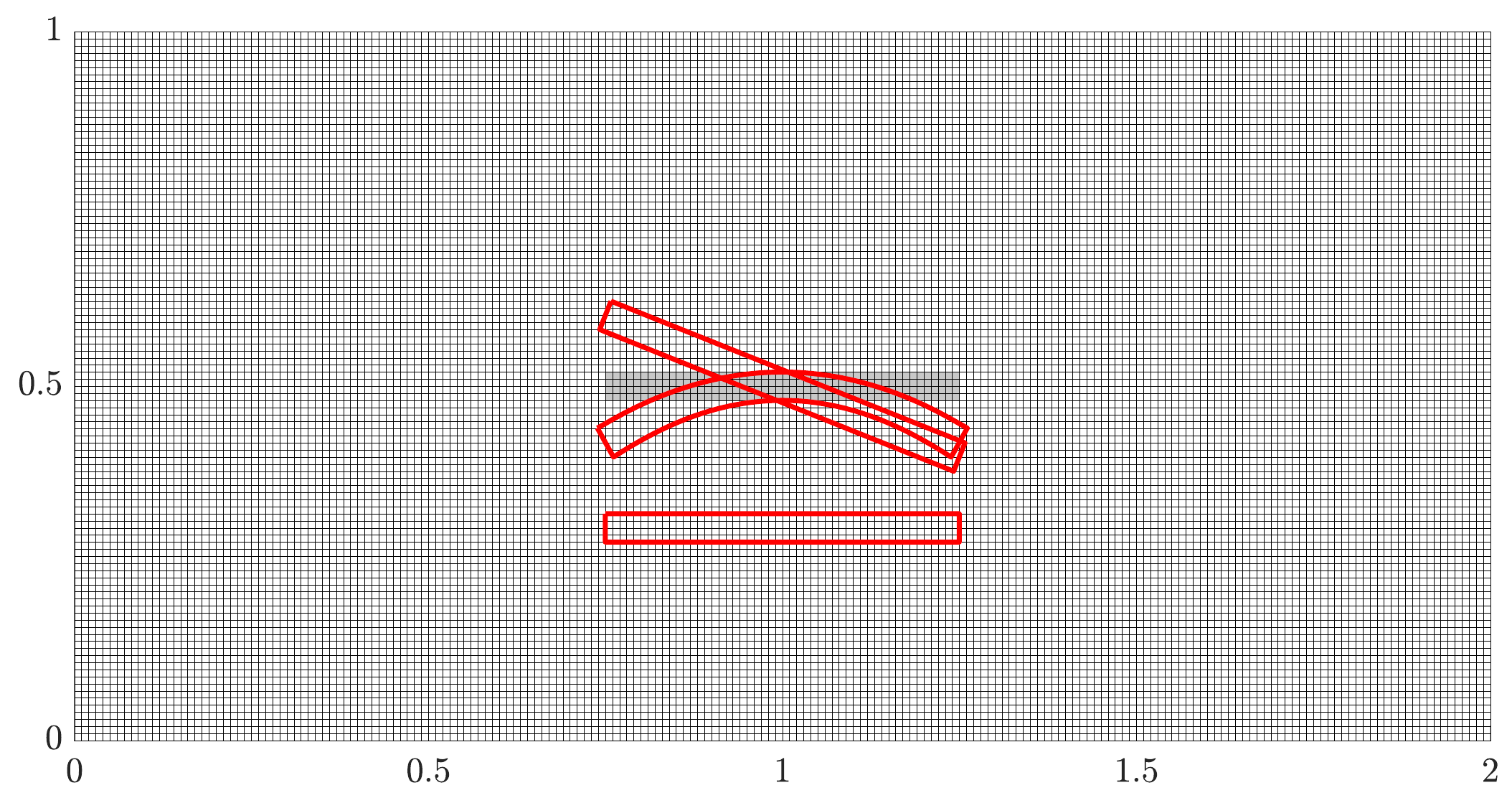}
	\caption{The original structured, quad mesh of the foil in a channel test problem. The original structure is shown in gray, while the deformed fluid-structure interface is shown in red. Dimensions are in meters unless otherwise noted.}
	\label{fig:foilinachanneltestprob}
\end{figure}

Before proceeding with our discussion of the different mesh deformation techniques, we state the test problems and mesh quality metrics to be used to demonstrate the functionality and tuning of each mesh deformation technique.  In the subsequent investigations of mesh quality, we focus mainly on the area surrounding the fluid-structure interface since this is where the main deformation occurs.

The \textbf{first test problem} is a standard benchmark problem in FSI literature; and is termed \textit{beam in a channel} (sometimes termed \textit{column or wall}) \citep[p.~604]{Yoon2010}. A cantilever beam is inserted inside a fluid channel such that its longitudinal axis is perpendicular to the flow in the initial state. The fluid flow exerts a force on the beam thus deflecting it. In this work, we hollowed out part of the base of the beam to increase the deflection at its base, and the deformed fluid-structure interface is in fact generated by running an FSI analysis. Figure \ref{fig:beaminachanneltestprob} shows the original structured, quad mesh and the deformed fluid-structure interface.

The \textbf{second test problem} is inspired by \citep[p.~61]{Stein2003}, and is termed \textit{foil in a channel} (since in the initial state, the longitudinal axis of the structure is parallel to the flow like an airfoil). In this test problem, three different types of deformations are applied separately; a rigid-body translation in the negative $y$ direction, a rigid-body rotation around the foil center of mass in the clockwise direction, and a prescribed U-shaped bending facing downwards. Figure \ref{fig:foilinachanneltestprob} shows the original structured quad mesh and the deformed fluid-structure interface for the three deformation modes overlaid in one figure. Unlike the first test problem, the deformation modes in the second test problem are generated by applying prescribed displacements to the foil, not as a result of an FSI analysis.

As for the \textbf{mesh quality metrics}, we utilize two metrics in this work; \textcolor{black}{\textit{skewness}} to describe the shape changes and \textcolor{black}{\textit{change in area}} to describe volume change. Skewness is a measure of the angular skewness of quad elements. \textcolor{black}{We adopt the definition in \citep[p.~677]{comsol56} where for} each finite element, skewness is calculated for each of its four corner angles $\theta$ as follows:
\begin{equation}
	1 - \textrm{max} \left( \frac{\theta - 90}{90}, \frac{90 - \theta}{90} \right)
	\label{eq:skewness}
\end{equation}
\noindent then the minimum is taken as the skewness measure for this element. This measure is also capable of detecting inverted elements for which it will give a negative value. Unlike triangular elements, the aspect ratio is not of critical value for quad elements. For instance, quad elements used in boundary layers have a large aspect ratio and is not considered a defect. \textcolor{black}{As for the second quality metric, the area change has been defined before for triangular elements in \cite{Johnson1996}. For the quad elements implements in this study, the area change is calculated simply as follows:}
\begin{equation}
	\frac{A_c}{A_0}
\end{equation}
\noindent where $A_c$ and $A_0$ are the current and original areas respectively. We pay attention to both the minimum and the maximum area changes.

\textcolor{black}{Even though boundary layers are significant in any work with fluid flow in contact with a no-slip boundary, a lot of work on TO of fluid flow problems still uses regular meshing. Two main reasons for this particular trend are:
	\begin{enumerate}
		\item Most of the works that address fluidic TO (or TOFSI) of high Reynolds numbers are still in the academic proof-of-concept stage and well beyond real-life applicability, hence the focus on resolving the boundary layers in the optimized designs is not a priority \textit{so far}.
		\item Most of the research studies so far is focused on applications where boundary layers shear stresses are not of great significance such as minimizing total energy dissipation in pure fluidic TO (\cite{Yoon2016, Yoon2020a}) or studying the effect of fluid flow on a structure such that pressure is more significant than shear stresses at the fluid-structure interface (\cite{Lundgaard2018, Yoon2017, Jenkins2016}).
	\end{enumerate}
In addition, since the fluid-structure interface is inherently implicit in TO problems (i.e. not known apriori), it would require severe mesh deformations at best or frequent remeshing at worst to resolve the boundary layers in the optimized designs.
}


\section{Legacy Mesh Deformation Techniques for TOFSI Problems}
In this section, we discuss and comment on some legacy mesh deformation techniques used in TOFSI. The representation of the fluid mesh as a continuous or discrete fictitious media (elastic or otherwise) is not new, it advanced hand in hand with the numerical modeling of problems that require a distinct fluid interface such as free surface flows and FSI.

\subsection{Discrete \textcolor{black}{Spring Analogy} Mesh Deformation Model}
\label{ssec:lineartorsionalsprings}
Historically, the deformation of the fluid mesh was used first for fluid flow problems with discontinuities (e.g., shock waves) or with multiple phases (e.g., free surface flow and flame propagation). \citet{Batina1990} used lineal\footnote{\textit{Lineal} is used to refer to straight springs as opposed to \textit{torsional}, while \textit{linear} is used to denote that the stiffness is linearly proportional to the displacement as opposed to \textit{non-linear}.} tension/compression springs along the finite element edges connecting the mesh nodes, where the spring stiffness in inversely proportional to the edge length. The spring stiffness can also be linked to a fluid flow criteria such as the Mach number. \citet{Palmerio1994} included a \say{pseudo-pressure} term to prevent local mesh size vanishing. Nonetheless, although this model prevents two mesh vertices from colliding, it doesn't prevent mesh vertices from crossing opposite edges in largely-deformed two-dimensional triangular meshes (i.e., creating inverted elements). \textcolor{black}{Later, \citet{Farhat1998} complemented the lineal springs model of \cite{Batina1990} with torsional springs at the vertices such that the torsional stiffness is calculated based on contributions from all the triangular elements sharing that particular vertex. The model by \cite{Farhat1998} is the focus of our work, and is henceforth referred to as the \say{the spring analogy model}.} This technique was later adapted to three-dimensional meshes in a clever indirect approach in \cite{Degand2002}. This method is implemented in the design optimization of aeroelastic systems \citep{Allen2002, Allen2004} and the topology optimization of aeroelastic structures \citep{Maute2002, Maute2004, Maute2006}\footnote{The work by \citet{Maute2002} is most probably the first to address topology optimization of fully-coupled fluid-structure interactions. Although the fluid-structure interface was excluded from the design space (i.e., dry topology optimization), this work was the first to employ three field sensitivity analysis in a topology optimization context.}. In all the above work, this technique is applied to unstructured, triangular fluid meshes.

\textcolor{black}{The spring analogy model is formulated as follows:
\begin{equation}
	\begin{aligned}
		(\mathbf{K}^{ij}_{lineal} + \mathbf{K}^{ijk}_{torsional}) \mathbf{x} & = \mathbf{0} \quad \mathrm{on} \quad \mathrm{\Omega} \backslash \mathrm{\Gamma}_{\mathrm{FSI}}, \\
		\mathbf{x} & = \mathbf{u} \quad \mathrm{on} \quad \mathrm{\Gamma}_{\mathrm{FSI}}.
	\end{aligned}
\end{equation}
\noindent where the stiffness matrices $\mathbf{K}^{ij}_{lineal}$ and $\mathbf{K}^{ijk}_{torsional}$ represent the contributions from the linear and torsional springs respectively, $\mathbf{x}$ is the mesh nodal displacements vector, $\mathrm{\Omega}$ is the mesh computational domain, $\mathrm{\Gamma_{\mathrm{FSI}}}$ is the fluid-structure interface, and $\mathbf{u}$ is the structural nodal displacements vector.}

\textcolor{black}{Each finite element edge is represented by a lineal spring that is assumed to be a generally oriented truss element of an angle $\alpha$ where the cross-sectional area and the material elastic modulus are assumed unity. Accordingly, the stiffness of each lineal spring is inversely proportional to its length $L$. The torsional springs part of the model is designed such that the torsional stiffness at each node depends on a combination of edge lengths $L_{ij}, \, L_{jk}, \, \textrm{and } L_{ki}$, and areas $A_{ijk}$ of all finite elements sharing this specific node as in Eq. \ref{eq:springequation}, where $L_{ij}$ is the edge length between nodes $i$ and $j$, $x_{ij}$ and $y_{ij}$ are the distance between nodes $i$ and $j$ along the $x$ and $y$ axes respectively, and $A_{ijk}$ is the area of the finite element sharing nodes $i$, $j$, and $k$ (cf. \cite{Farhat1998} for the detailed derivation of Eq.\ref{eq:springequation}).}

Next, we discuss in detail some important points about this model.

\begin{strip}
	\begin{equation}
		\begin{aligned}
			\mathbf{K}^{ij}_{lineal} & = \frac{1}{L_{ij}} \begin{bmatrix}
				\cos[2](\alpha) & \sin(\alpha) \cos(\alpha) & - \cos[2](\alpha) & - \sin(\alpha) \cos(\alpha) \\
				\sin(\alpha) \cos(\alpha) & \sin[2](\alpha) & - \sin(\alpha) \cos(\alpha) & - \sin[2](\alpha) \\
				- \cos[2](\alpha) & - \sin(\alpha) \cos(\alpha) & \cos[2](\alpha) & \sin(\alpha) \cos(\alpha) \\
				- \sin(\alpha) \cos(\alpha) & - \sin[2](\alpha) & \sin(\alpha) \cos(\alpha) & \sin[2](\alpha)
			\end{bmatrix}_{ij}, \\
			\mathbf{K}^{ijk}_{torsional} & = {\mathbf{R}^{ijk}}^T \, \mathbf{C}^{ijk} \, \mathbf{R}^{ijk}, \\
			\mathbf{R}^{ijk} & = \begin{bmatrix}
				b_{ik} - b_{ij} & a_{ij} - a_{ik} & b_{ij} & - a_{ij} & - b_{ik} & a_{ik} \\
				- b_{ji} & a_{ji} & b_{ji} - b_{jk} & a_{jk} - a_{ji} & b_{jk} & - a_{jk} \\
				b_{ki} & - a_{ki} & - b_{kj} & a_{kj} & b_{kj} - b_{ki} & a_{ki} - a_{kj}
			\end{bmatrix},\\
			a_{ij} & = \frac{x_{ij}}{L^2_{ij}}, \\
			b_{ij} & = \frac{y_{ij}}{L^2_{ij}}, \\
			\mathbf{C}^{ijk} & = \begin{bmatrix}
				C^{ijk}_i & 0 & 0 \\
				0 & C^{ijk}_j & 0 \\
				0 & 0 & C^{ijk}_k
			\end{bmatrix}, \\
			C^{ijk}_i & = \frac{L^2_{ij} L^2_{ik}}{4 A_{ijk}^2}. \label{eq:springequation} 
		\end{aligned}
	\end{equation}
\end{strip}

\subsubsection{Non-Linearity and Computational Costs}
\label{secss:compcosts}
In the derivation of the matrix $\mathbf{R}^{ijk}$ in Eq. \ref{eq:springequation}, an approximation is use\textcolor{black}{,} $\sin(\Delta \theta) \approx \Delta \theta$\textcolor{black}{,} to obtain a linear relation between the angles of each triangle and the displacement of its vertices \citep[p.~234]{Farhat1998}. This approximation necessitates solving the mesh deformation equations in a number of smaller steps even if there are no other constraints on the step size as is the case with quasi-static problems. This results in a big disadvantage for this model; that is its computational cost. Potentially, obtaining an exact linear relation between the angles and the displacements - if any exists - could be a game changer for this model. But again, probably its non-linearity is one of the reasons why it behaves well. From our numerical experiments on the test problems, as small as 5 time steps are enough to prevent any overlapping or degenerated elements, but somewhere in the vicinity of 30 time steps are needed to enhance the minimum element quality depending on the problem.

\subsubsection{Application to Quadrilateral Elements}

\begin{table*}
	\small
	\centering
	\captionsetup{width=0.88\textwidth}
	\caption{Effect of diagonal selection on the minimum skewness using the spring analogy model.}
	\label{tab:triangulation}
	\begin{tabular}{@{}lcccc@{}} \toprule
		& \multicolumn{4}{c}{Minimum Skewness} \\
		& \multicolumn{4}{c}{\textcolor{blue}{Minimum Area Change}} \\ 
		& \multicolumn{4}{c}{\textcolor{red}{Maximum Area Change}} \\ 
		\cmidrule{2-5}
		Test Problem & Diagonal 1-3 & Diagonal 2-4 & Both Diagonals & Selective Diagonals \\ \midrule
		Beam in a Channel & 0.065 & 0.283 & 0.241 & 0.327 \\
		& \textcolor{blue}{0.247} & \textcolor{blue}{0.280} & \textcolor{blue}{0.264} & \textcolor{blue}{0.239} \\
		& \textcolor{red}{2.740} & \textcolor{red}{2.764} & \textcolor{red}{2.733} & \textcolor{red}{2.890} \\
		Foil in a Channel - Translation & 0.092 & 0.092 & 0.243 & 0.338 \\
		& \textcolor{blue}{0.226} & \textcolor{blue}{0.226} & \textcolor{blue}{0.232} & \textcolor{blue}{0.227} \\
		& \textcolor{red}{3.046} & \textcolor{red}{3.046} & \textcolor{red}{3.009} & \textcolor{red}{3.247} \\
		Foil in a Channel - Rotation & 0.320 & -0.022 & 0.211 & 0.295 \\
		& \textcolor{blue}{0.300} & \textcolor{blue}{0.274} & \textcolor{blue}{0.288} & \textcolor{blue}{0.274} \\
		& \textcolor{red}{2.755} & \textcolor{red}{2.714} & \textcolor{red}{2.708} & \textcolor{red}{2.970} \\
		Foil in a Channel - Bending & 0.008 & 0.008 & 0.245 & 0.313 \\
		& \textcolor{blue}{0.339} & \textcolor{blue}{0.339} & \textcolor{blue}{0.354} & \textcolor{blue}{0.321} \\
		& \textcolor{red}{2.406} & \textcolor{red}{2.406} & \textcolor{red}{2.356} & \textcolor{red}{2.472} \\	\bottomrule
	\end{tabular}
\end{table*}

The spring analogy model in \cite{Farhat1998} is derived for tri elements, hence it is not suitable for application directly to quad elements in its current form. Thus it's imperative that all quad elements are triangulated before applying this model. For every quad element, there are three types of triangulations; along either diagonals separately or along both diagonals at the same time. For instance, for a quad element with nodes ordered counter clockwise, the possible triangulations are: connect nodes 1 and 3, connect nodes 2 and 4, or connect both. Given that this is a structured mesh, it can be easily arranged that node numbering for all elements start from the same position (bottom left and going counterclockwise in this work) and controlling the triangulation is easily performed from the nodal numbers directly. Alternatively, it could be beneficial to pick a suitable diagonal for each element individually. This could be easily performed by solving two time steps with opposing triangulations then picking the suitable diagonal for each element based on which produces the higher quality (minimum skewness in this work). This technique is henceforth termed \textit{selective diagonals}. Note that even though we work with tri elements in the spring analogy model, the skewness measure is still applied to the quad\textcolor{black}{,} not tri\textcolor{black}{,} elements. 

Table \ref{tab:triangulation} shows mesh quality metrics for the test problems using all possible combinations of diagonals. Selective diagonals produce better skewness in all test problems except the foil in a channel under rotation. This particular test problem has the same deformation mode for all elements in the vicinity of the fluid-structure interface and hence a particular diagonal is more suitable in this case. Nonetheless, the selective diagonals technique's skewness is only 8\% worse. The reason why selective diagonals didn't work best in this case is probably due to the dependence between adjacent elements; meaning that altering the diagonal in one element might affect the optimum diagonal orientation in an adjacent element. However, selecting the absolute optimum diagonal for each element would require setting up an iterative optimization problem which is a tedious task and we conjecture that the benefit in minimum quality increase is not worth it. Even though a quad element is inverted in the foil in a channel under rotation problem when using diagonal 2-4, it is no fault of the spring analogy model which successfully resisted inverting/overlapping any tri elements. The foil in a channel problem under translation and under bending is mirror symmetric around a vertical axis, hence either diagonal individually gives the same low quality. The beam in a channel test problem performed best for the selective diagonals technique. This is probably due to the fact that it contains a mixture of deformation modes unlike the other test cases. Considering that only the skewness measure was used in selecting the optimal diagonals, selective diagonals doesn't work well w.r.t. the area change mesh metric. A weighted combination of skewness and area change could potentially be used in selecting the optimal diagonal for each finite element but this compromise is inherent in the nature of mesh deformation and is better left to the judgment of the user.

\begin{figure*}
	\centering
	\subfloat[Diagonal 1-3. \label{figs:diag1}]{\includegraphics[width=0.245\linewidth]{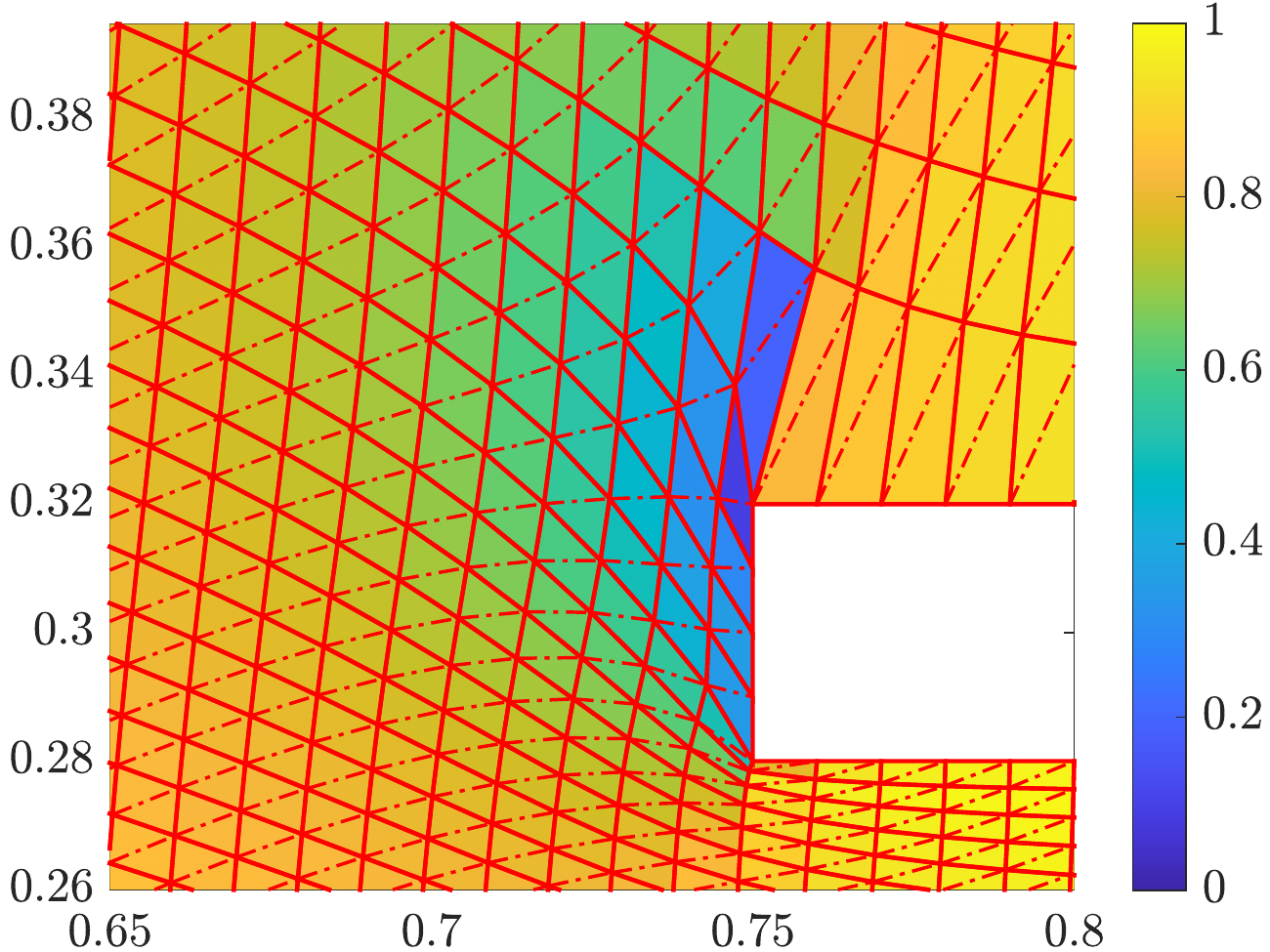}}
	\subfloat[Diagonal 2-4. \label{figs:diag2}]{\includegraphics[width=0.245\linewidth]{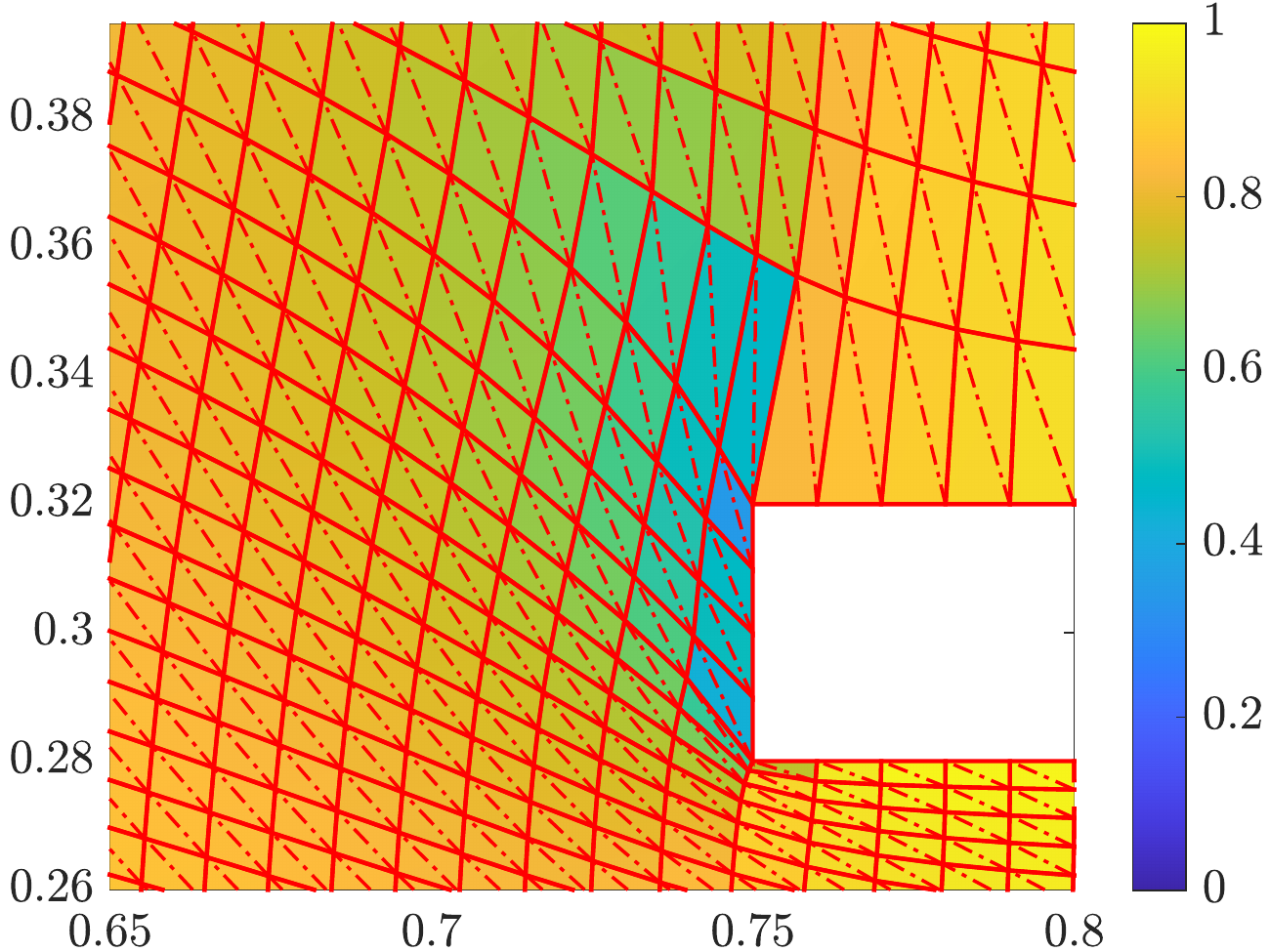}}
	\subfloat[Both Diagonals. \label{figs:diag3}]{\includegraphics[width=0.245\linewidth]{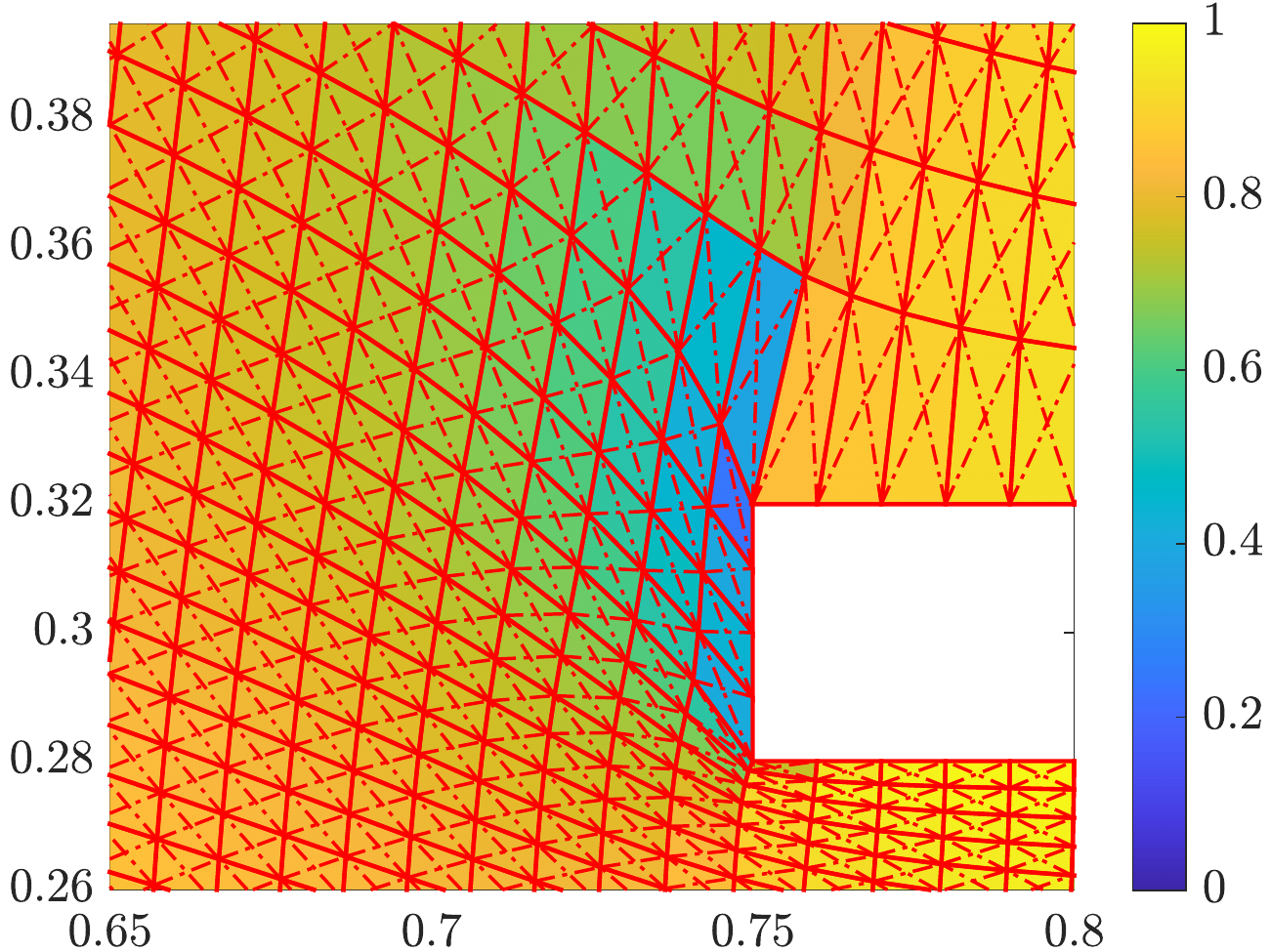}}
	\subfloat[Selective Diagonals. \label{figs:diag4}]{\includegraphics[width=0.245\linewidth]{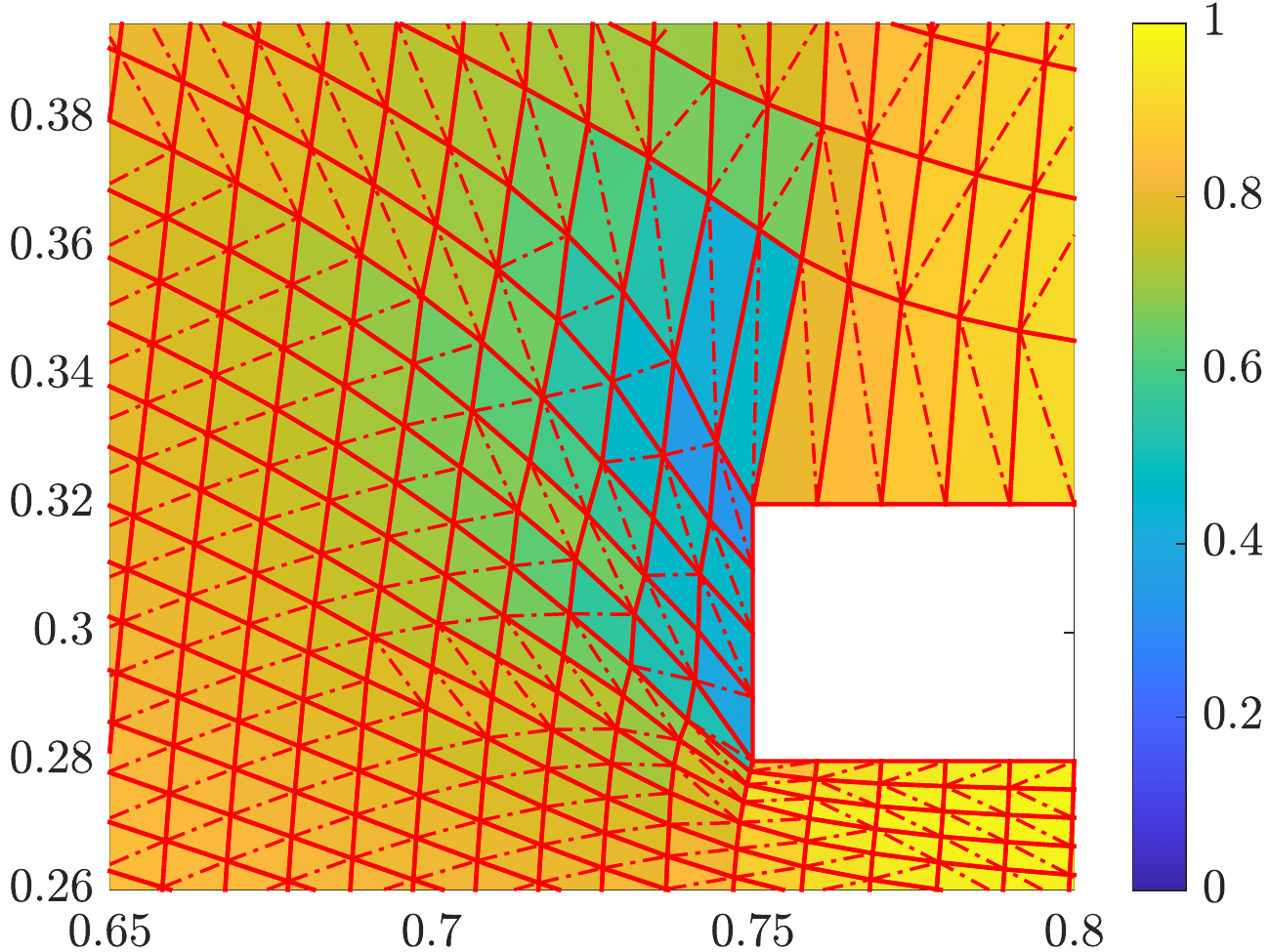}}
	\caption{Element skewness in the foil in a channel test problem under a downward translation mode. Element edges are shown in solid red, while the diagonals are shown in dash-dot red.}
	\label{fig:quad_triangulation}
\end{figure*} 

Figure \ref{fig:quad_triangulation} shows element skewness in the foil in a channel test problem under translation with a focus on elements located to the left of the foil, which are subjected to almost pure shear. Figure \ref{figs:diag1} shows a triangulation along nodes 1 \& 3, it doesn't align well with this shear stress mode and the elements shown are heavily skewed. On the other hand, Fig. \ref{figs:diag2} shows a better triangulation along nodes 2 \& 4 and the elements shown support the shear stresses more effectively. Figure \ref{figs:diag3} shows both triangulations overlaid on top of each other and the result is somewhere in between Fig. \ref{figs:diag1} and Fig. \ref{figs:diag2}. Of course, triangulation along both diagonals simultaneously comes at almost double the computational cost of triangulating along either diagonal separately. Figure \ref{figs:diag4} shows the result of using selective diagonals for each element, which is better than all the other cases.


A few worthy remarks about the two initial test runs required for the selective diagonals technique are in order. \textbf{(i)} The additional cost of these two additional time steps would be negligible considering the total number of time steps needed as is discussed in \ref{secss:compcosts}. \textbf{(ii)} Since in this work the starting nodes for all elements are in the same position (bottom left), the same diagonals are used for all elements in each test run. However, if this is not the case, the result may be different, and this case was not tested here. \textbf{(iii)} The percentage of prescribed deformation used to solve the two initial test runs might have an effect on the selected diagonal for each element. It should neither be too small that it isn't representative of the real deformation modes, nor too large that the approximate solver is inaccurate.


\subsubsection{Dominance of the Lineal vs. Torsional Stiffness}
\begin{figure*}
	\centering
	\subfloat[GSC 0.001.]{\includegraphics[width=0.245\linewidth]{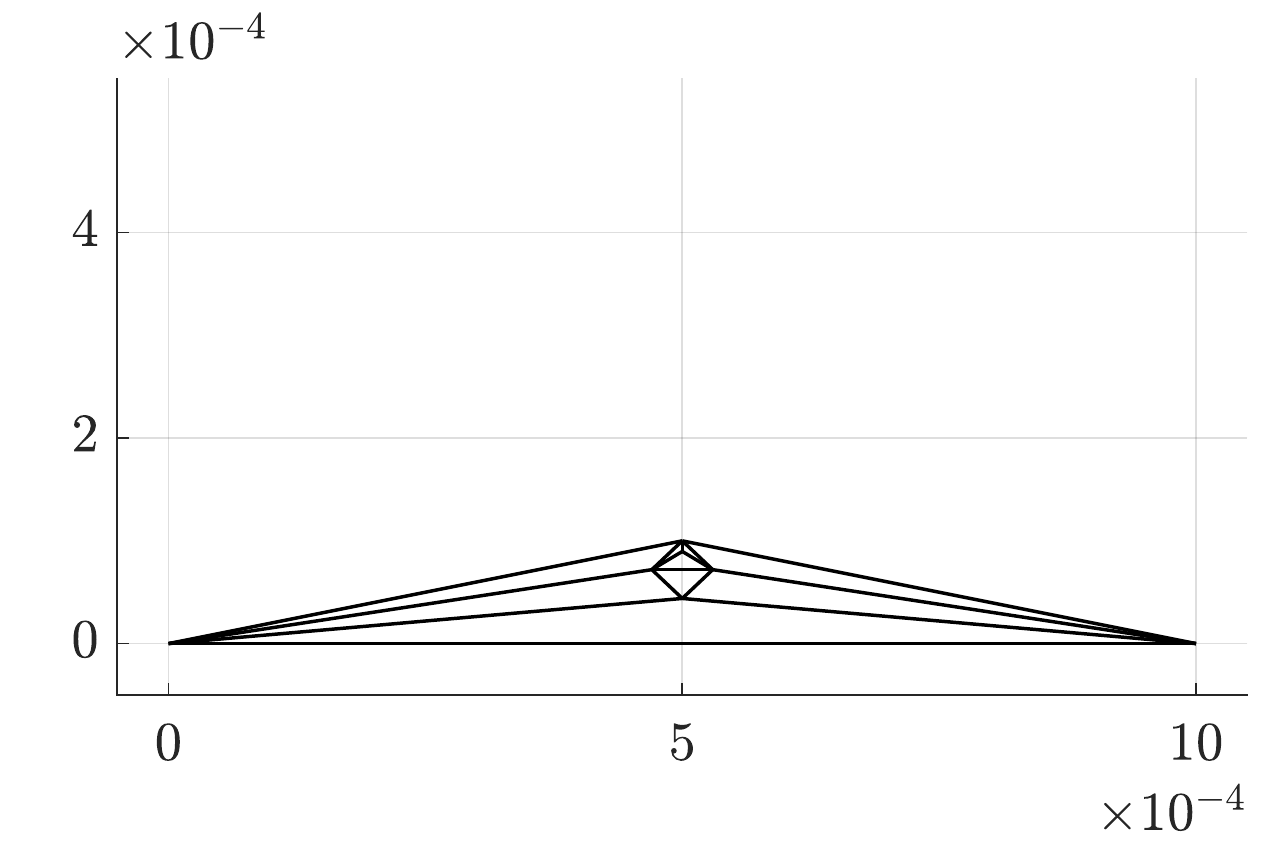}}
	\subfloat[GSC 1.]{\includegraphics[width=0.245\linewidth]{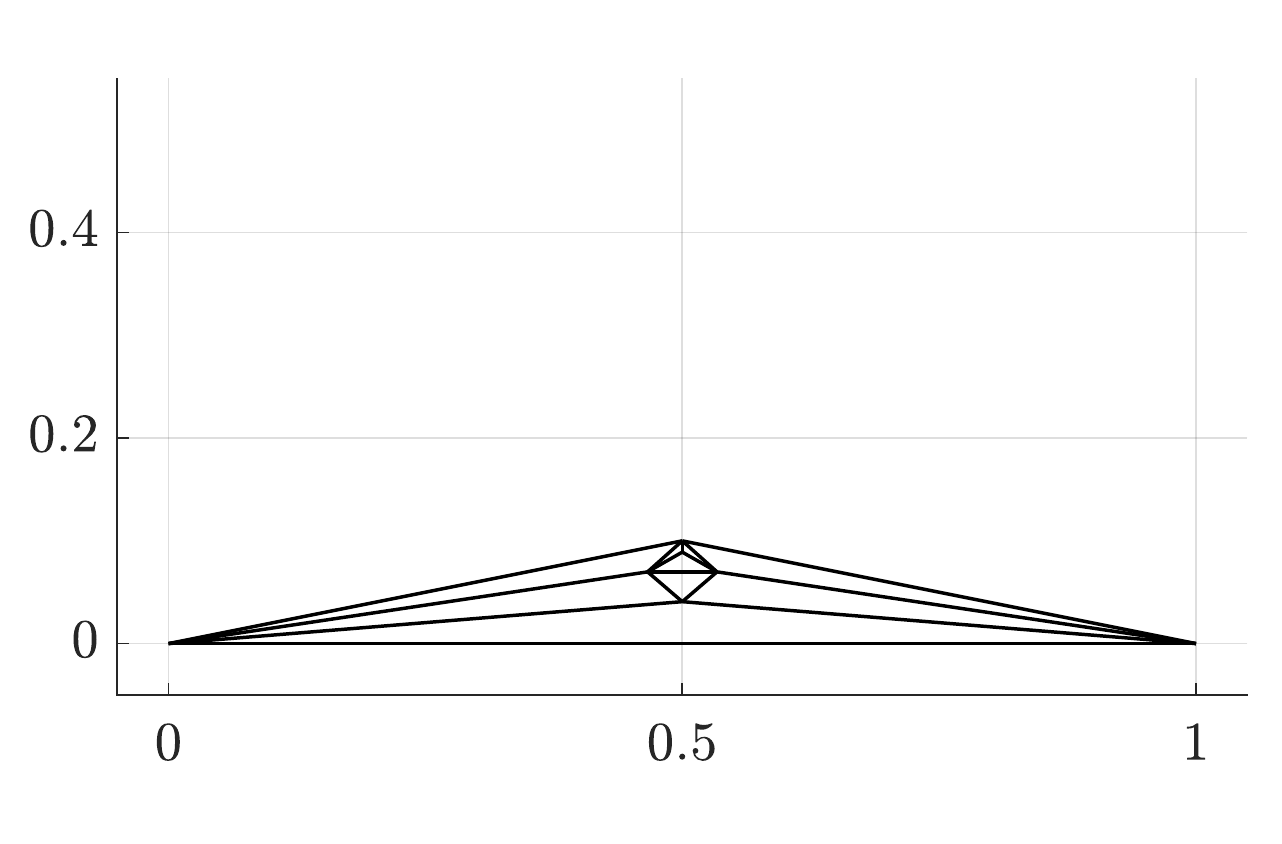}}
	\subfloat[GSC 100.]{\includegraphics[width=0.245\linewidth]{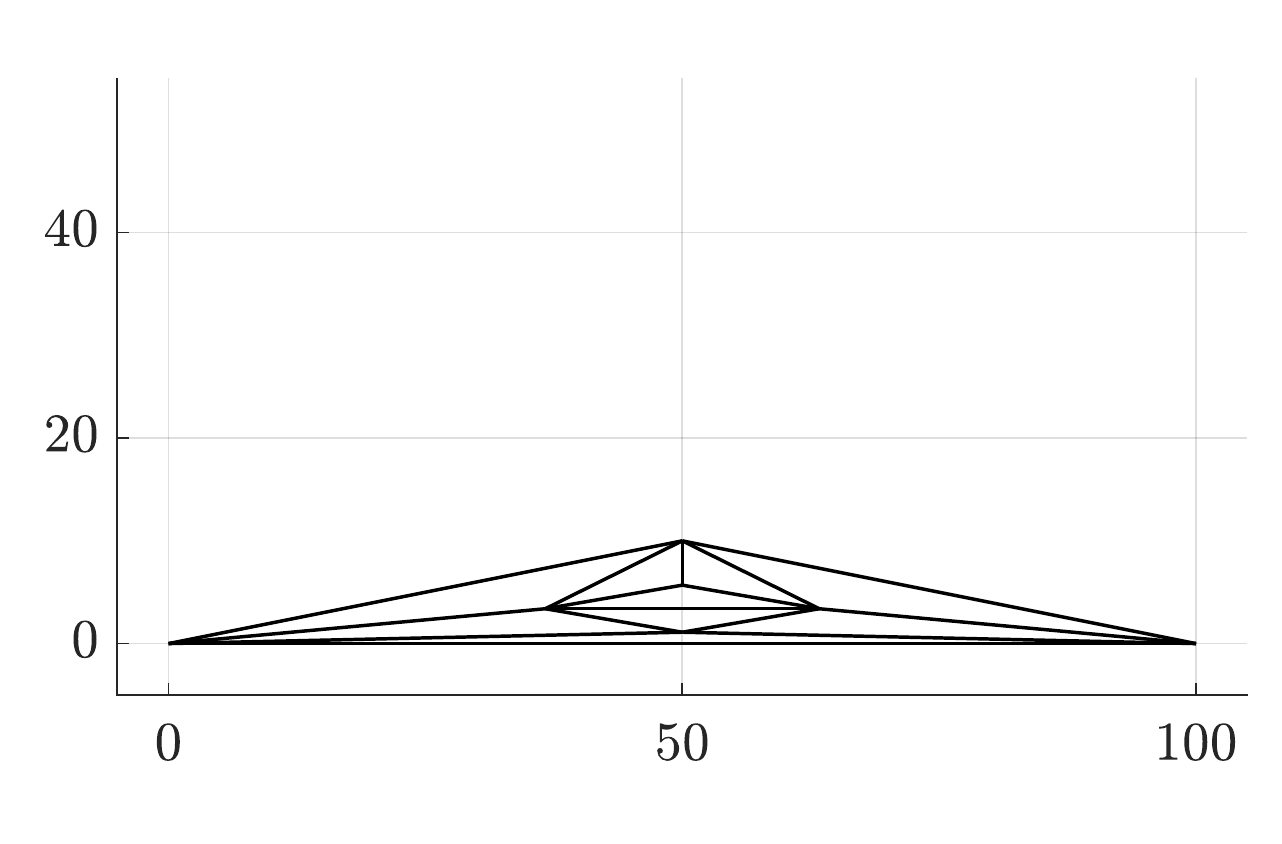}}
	\subfloat[GSC 1000.]{\includegraphics[width=0.245\linewidth]{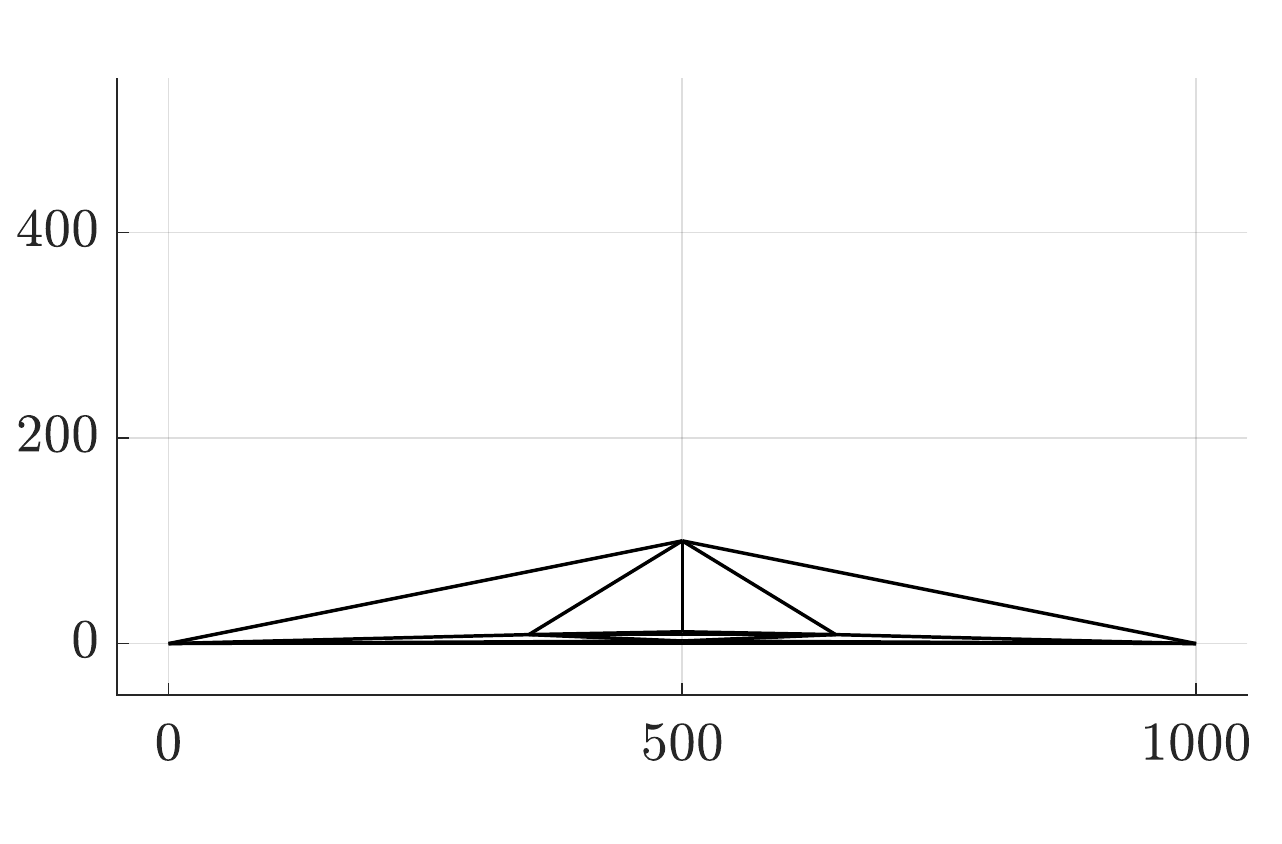}} \\
	\subfloat[TSC 1000.]{\includegraphics[width=0.245\linewidth]{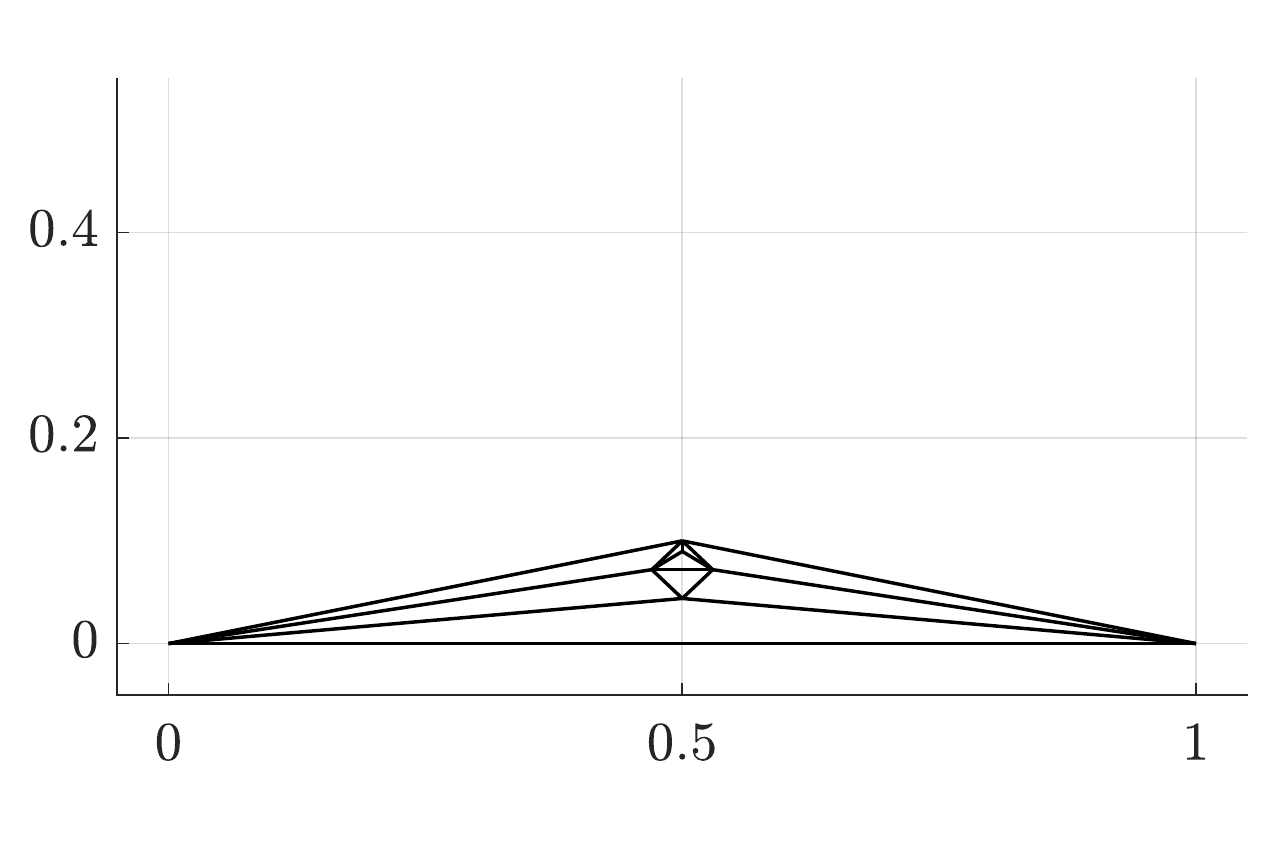}}
	\subfloat[TSC 1.]{\includegraphics[width=0.245\linewidth]{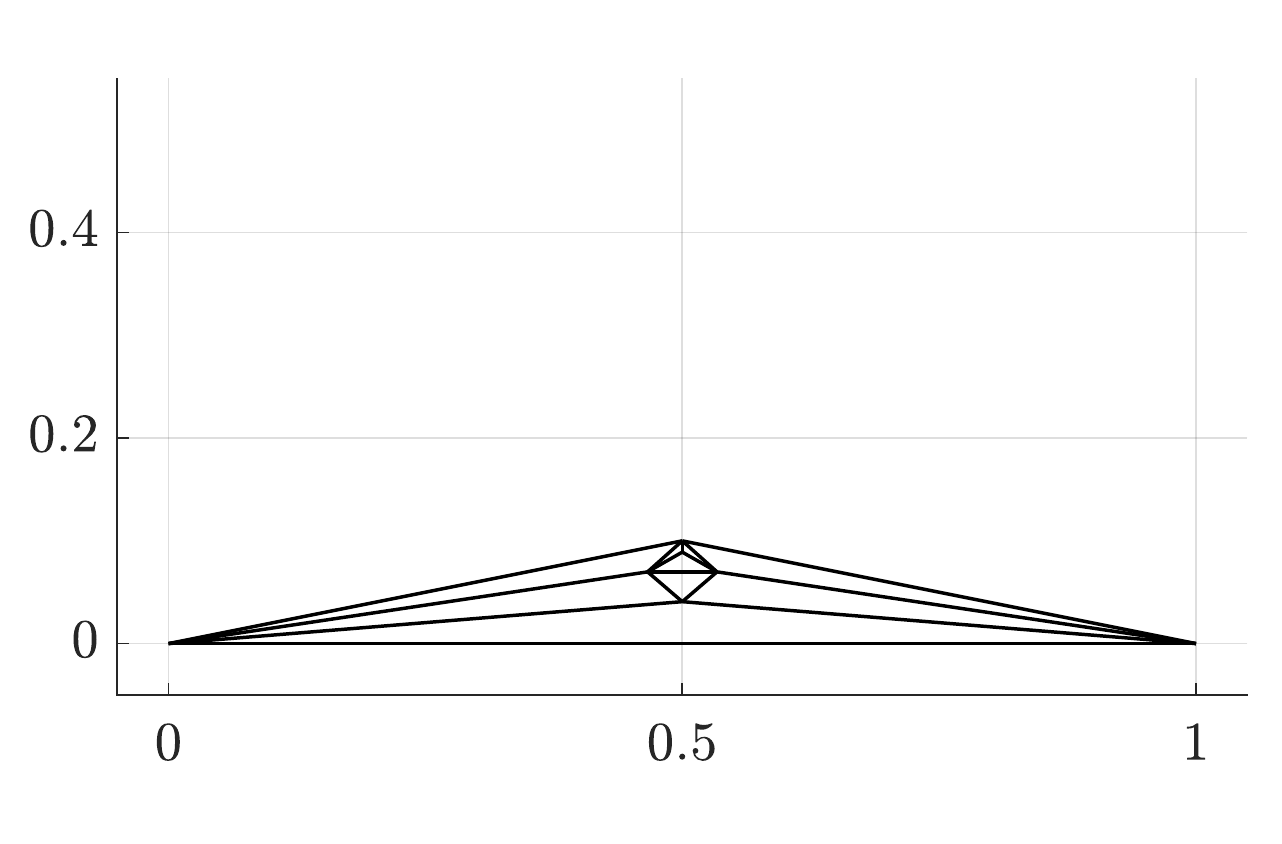}}
	\subfloat[TSC 0.01.]{\includegraphics[width=0.245\linewidth]{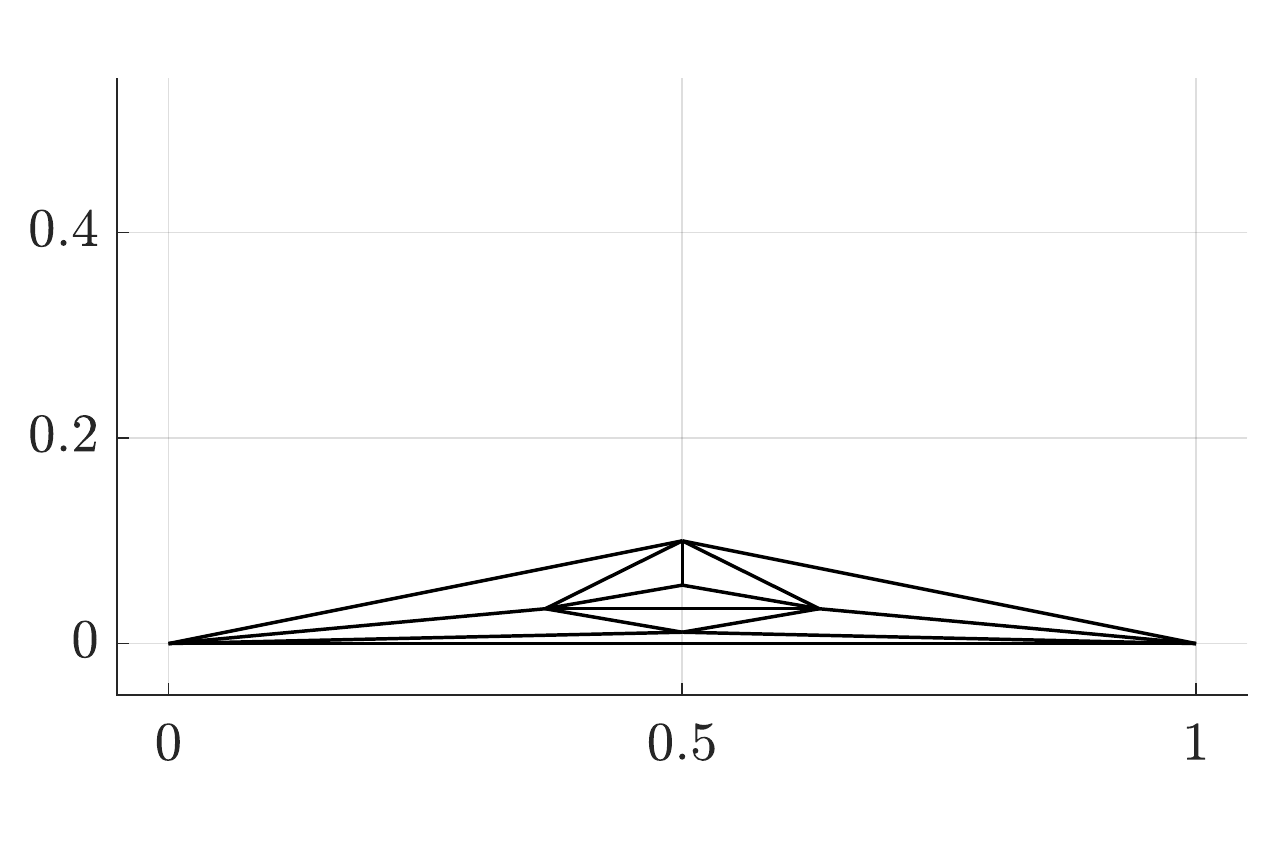}}
	\subfloat[TSC 0.001.]{\includegraphics[width=0.245\linewidth]{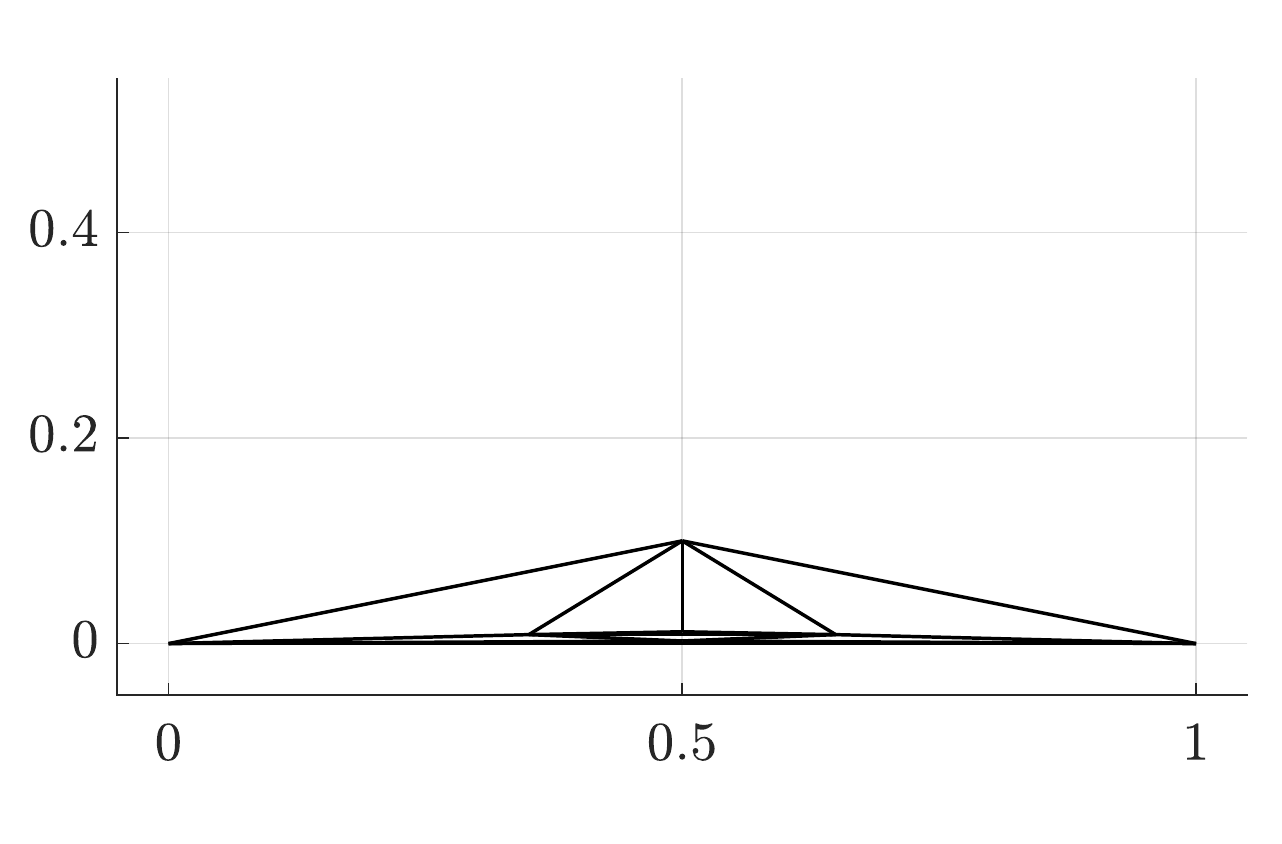}}
	\caption{Effect of the geometrical (GSC) and torsional (TSC) stiffness scaling factors on the behavior of the spring analogy model.}
	\label{fig:farhatexample}
\end{figure*}

Given the formulation of the spring analogy model in Eq. \ref{eq:springequation}, changing the geometric scale of the problem doesn't affect the behavior of the lineal springs part of the model. However,  the torsional springs part of the model - with its dependence on edge lengths and areas raised to different powers - is affected. That is changing the geometric size of the problem alters the relative contribution/dominance of the lineal vs. torsional stiffness of the model. To further elaborate on this issue, consider the original example mentioned in \citep[p.~239]{Farhat1998}, a large triangular shape consisting of 9 smaller triangles. The two bottom nodes are fixed, while a prescribed displacement in the negative $y$ direction is applied to the top node. In the top row of Fig. \ref{fig:farhatexample}, we alter the geometric scale by multiplying the nodal coordinates and the prescribed displacements by a factor (i.e., GSC). It can be observed that a high enough geometrical scale would cause inverted elements. The only reason inverted elements occur is if the torsional springs part of the model is undermined. In the bottom row of the same figure, we alter the scale of the torsional stiffness by multiplying $\mathbf{C}_{ijk}$ by a factor (i.e., TSC). It turns out that the geometrical scale is inversely proportional to the torsional stiffness scale due to the formulation of the torsional part in Eq. \ref{eq:springequation}. As expected, the lineal stiffness tends to resist changing the length of the finite element edges (left to right in Fig. \ref{fig:farhatexample}), while the torsional stiffness tends to resist the rotation of the element edges and hence resist the skewing/collapsing of the triangular elements (right to left in Fig. \ref{fig:farhatexample}). From a quad element perspective, this means that the lineal stiffness attempts to counteract volume changes while the torsional stiffness attempts to counteract shape changes. In most cases, shape changes are more critical than volume changes. Fortunately, the diminishing of the torsional stiffness only occurs at relatively large geometrical scales that are not common in real life applications.

\subsubsection{Layered Selective Stiffening at the Fluid-Structure Interface}
\label{secss:selective_stiffening_sprng}
Selective treatment of individual mesh deformation elements have been applied before in various forms to linear and non-linear elasticity models \citep{Tezduyar1992computation, Stein2003, Stein2004, Shamanskiy2021}. In this work, we apply \textit{layered selective stiffening} to the spring analogy model. We mainly stiffen \textcolor{black}{consecutive layers of} elements in the vicinity of the fluid-structure interface (i.e., the region where the flow solution is most critical). This simply shifts the distorted elements away from the fluid-structure interface where the bulk of mesh distortion is located.

\textcolor{black}{Before applying the layered selective stiffening, we need to sort mesh elements in consecutive layers starting from the fluid-structure interface, which is detailed as follows. First, the mesh nodes located on the fluid-structure interface are identified based on the overlapping between the solid and the fluid meshes. Second, the first layer of finite elements closest to the fluid-structure interface are identified through searching the solid mesh for elements that share nodes with the fluid-structure interface. Third, the second layer of finite elements are identified through searching the solid mesh for elements that share nodes with the first layer. The same procedure is followed for identifying further layers if needed. A worthy note specific to the spring analogy model is that the elements to be stiffened are selected as quads, then they are triangulated after.}

\begin{figure*}
	\centering
	\subfloat[\label{figs:spring_sel_stiff_one_layer_skew}]{\includegraphics[width=0.33\linewidth]{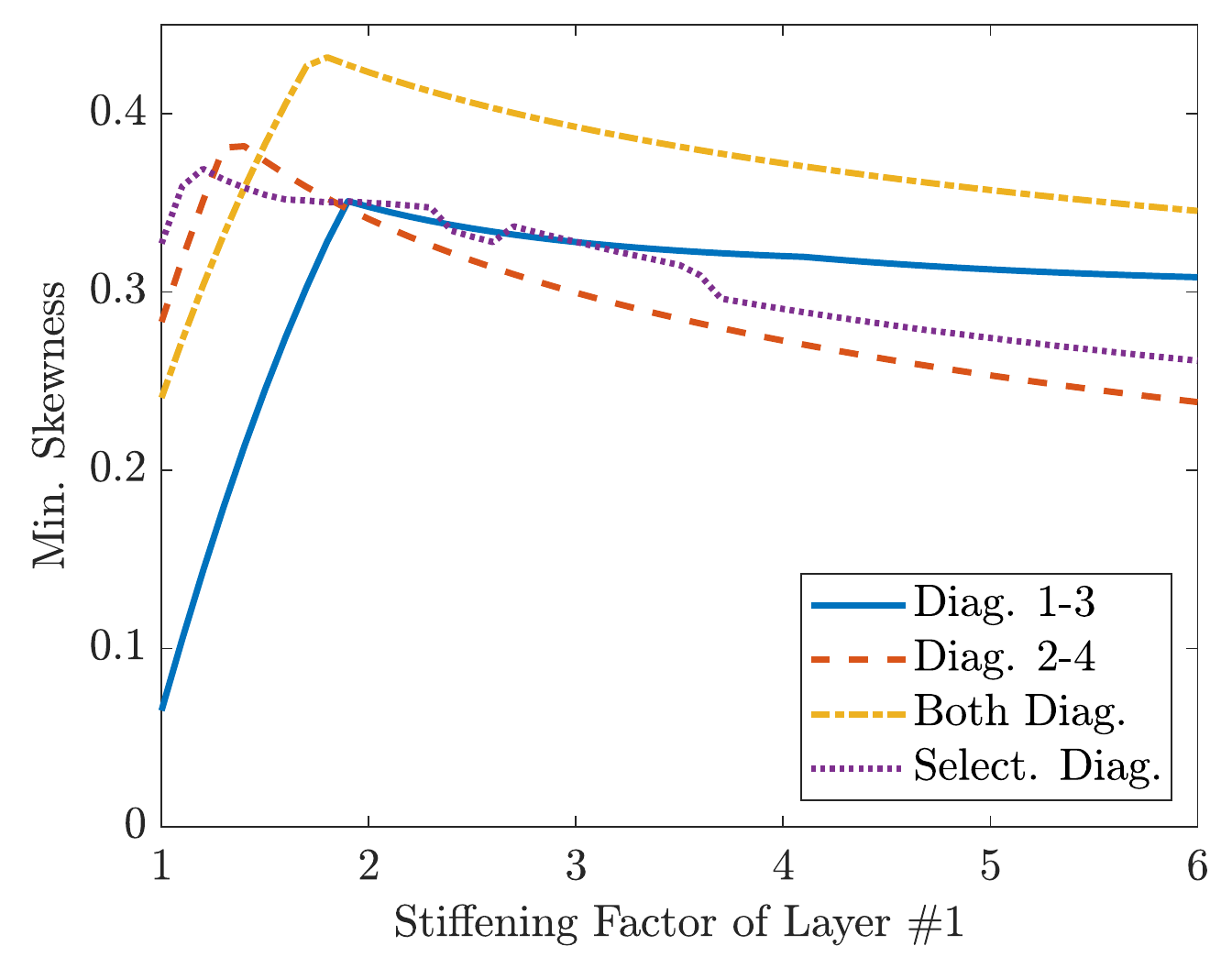}}
	\subfloat[\label{figs:spring_sel_stiff_one_layer_minar}]{\includegraphics[width=0.33\linewidth]{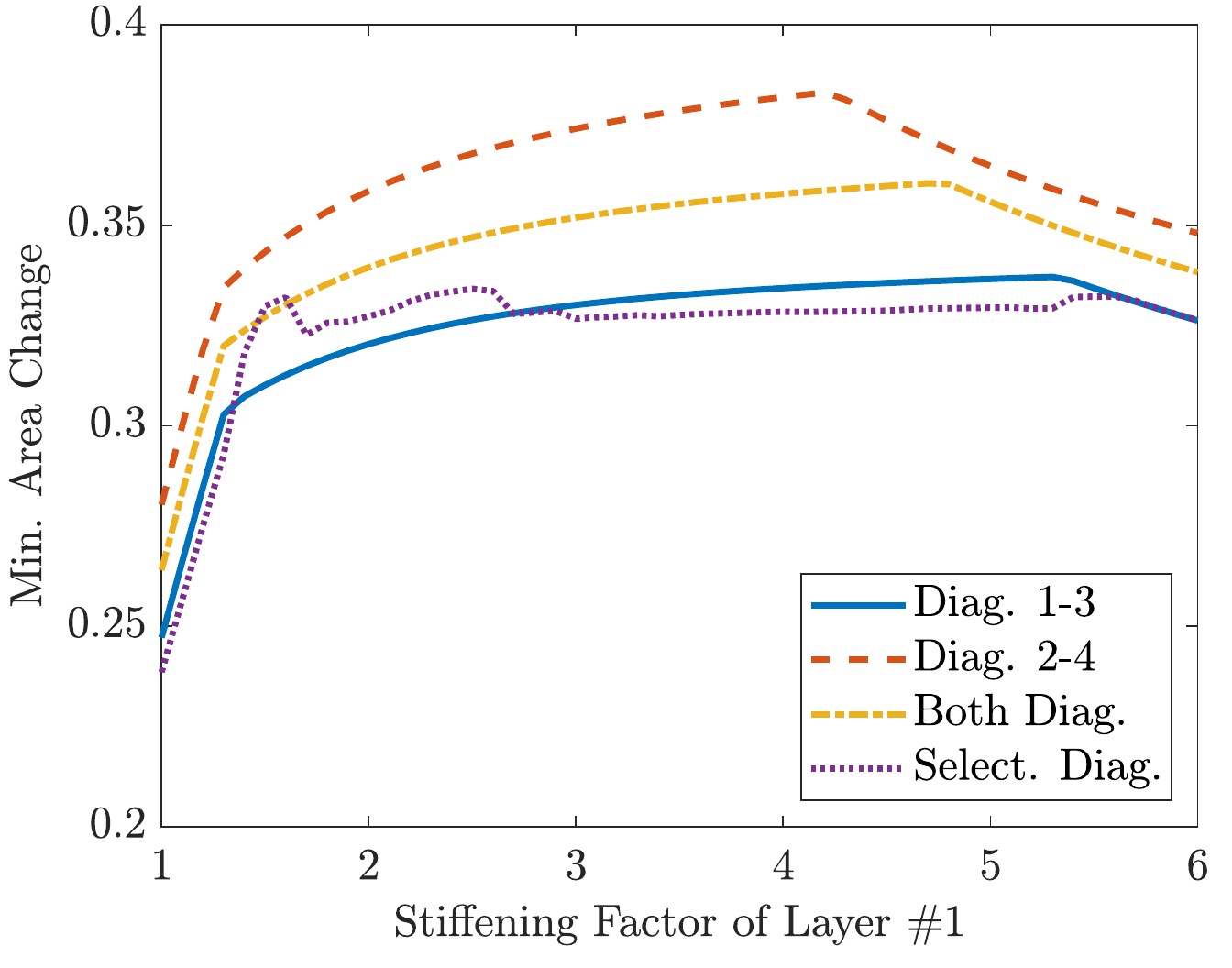}}
	\subfloat[\label{figs:spring_sel_stiff_one_layer_maxar}]{\includegraphics[width=0.33\linewidth]{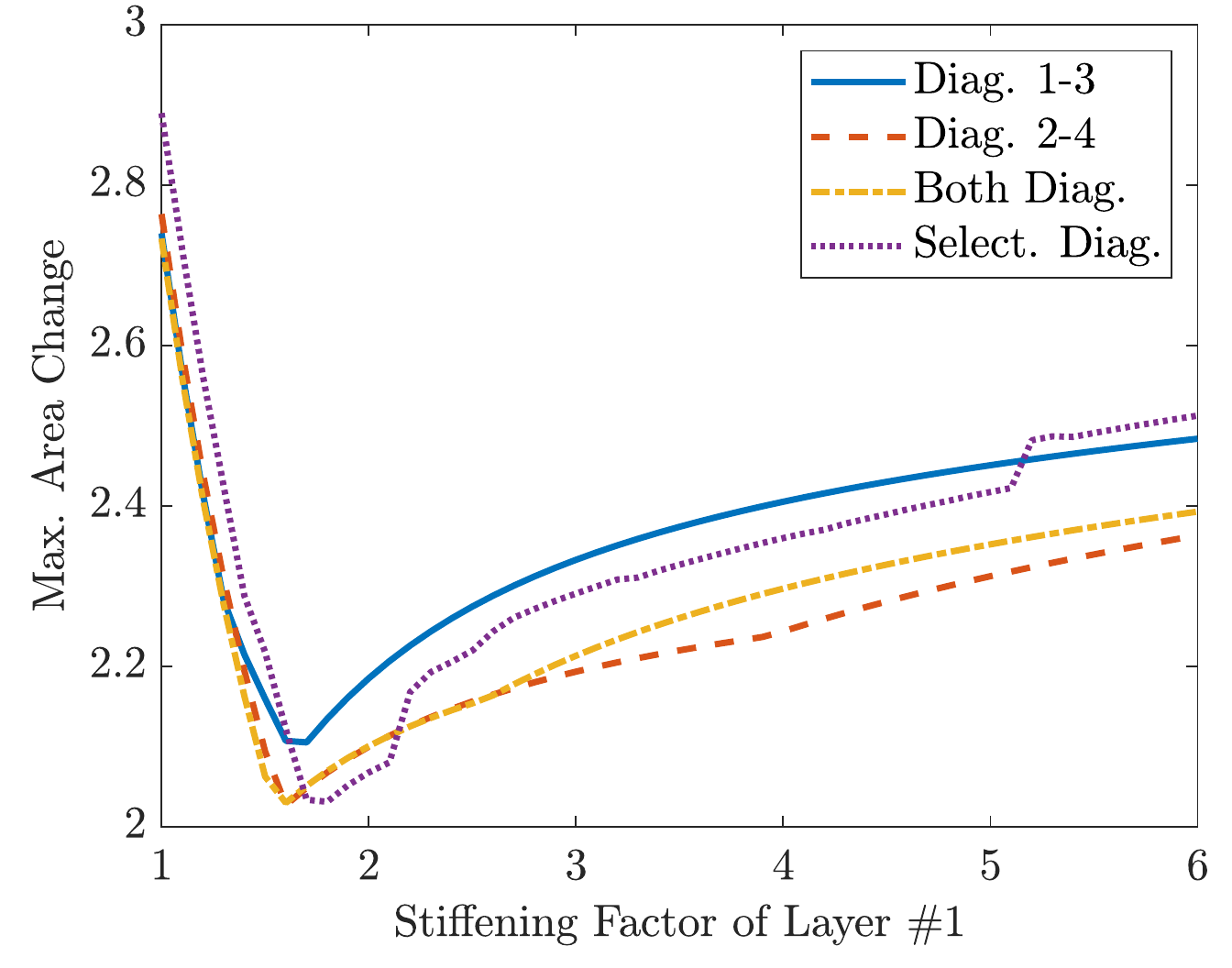}}
	\caption{Effect of selective stiffening of one layer on the mesh quality of the beam in a channel test problem using the spring analogy model.}
	\label{fig:spring_sel_stiff_one_layer}
\end{figure*}

We implement \textit{layered selective stiffening} by multiplying the local stiffness matrix (both lineal and torsional) by a stiffening factor greater than unity. Figure \ref{fig:spring_sel_stiff_one_layer} presents the effect of \textbf{selective stiffening of one layer} on the mesh quality of the spring analogy model for the beam in a channel test problem. Discontinuities in the graphs of the selective diagonals technique are due to diagonals switching as a result of changing the stiffening factor. Although for the case without stiffening, the selective diagonals technique performs best w.r.t. skewness, both diagonals together perform even better with stiffening. In Fig. \ref{figs:spring_sel_stiff_one_layer_skew}, the best skewness is achieved around a 1.8 stiffening factor. Similar minimum skewness results are obtained for the other test problems (Table \ref{tab:onelayerstiff_spranalg}), except the foil in a channel under rotation where diagonal 1-3 still performs best with stiffening. While the maximum area change  (Fig. \ref{figs:spring_sel_stiff_one_layer_maxar}) shows a minimum (i.e., an optimum), the minimum area change  (Fig. \ref{figs:spring_sel_stiff_one_layer_minar}) shows a monotonous increase with the stiffening factor, although at a diminishing rate. It's clear that there has to be a compromise between mesh quality metrics based on which is more critical for each test problem.

\begin{table*}
	\small
	\centering
	\captionsetup{width=0.85\textwidth}
	\caption{Effect of one layer selective stiffening on the minimum skewness using the spring analogy model.}
	\label{tab:onelayerstiff_spranalg}
	\begin{tabular}{@{}lccll@{}} \toprule
		Test Problem & Minimum Skewness & Stiffening Factor & Type of Diagonals \\ \midrule
		Beam in a Channel & 0.432 & 1.80 & Both Diagonals \\
		Foil in a Channel - Translation & 0.420 & 1.70 & Both Diagonals \\
		Foil in a Channel - Rotation & 0.478 & 4.10 & Diagonal 1-3 \\
		Foil in a Channel - Bending & 0.435 & 3.00 & Both Diagonals \\ \bottomrule
	\end{tabular}
\end{table*}

\begin{figure*}
	\centering
	\subfloat[\label{figs:SrngAnlg_sel_stiff_two_layers_big_step}]{\includegraphics[width=0.495\linewidth]{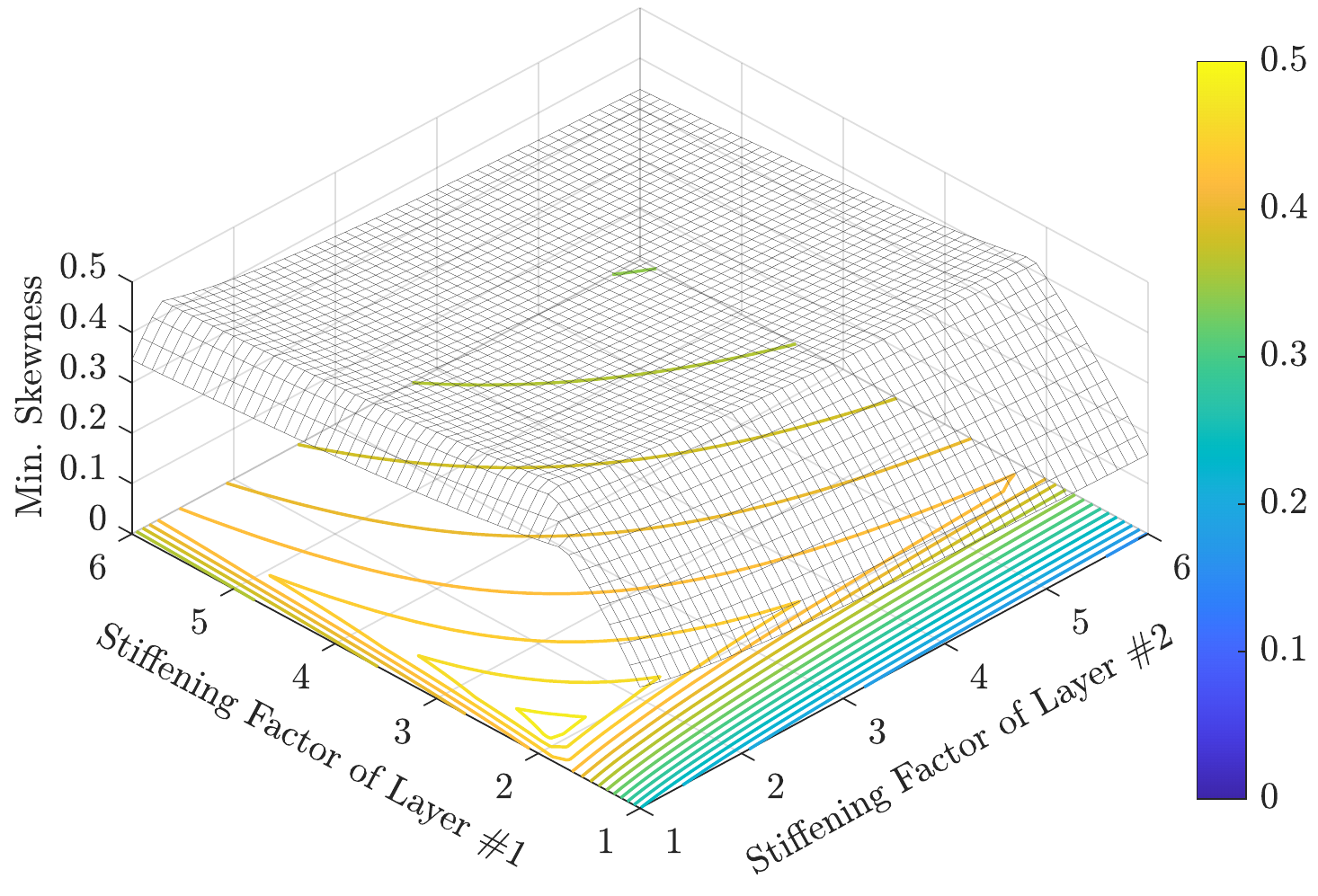}}
	\subfloat[\label{figs:SrngAnlg_sel_stiff_two_layers_small_step}]{\includegraphics[width=0.495\linewidth]{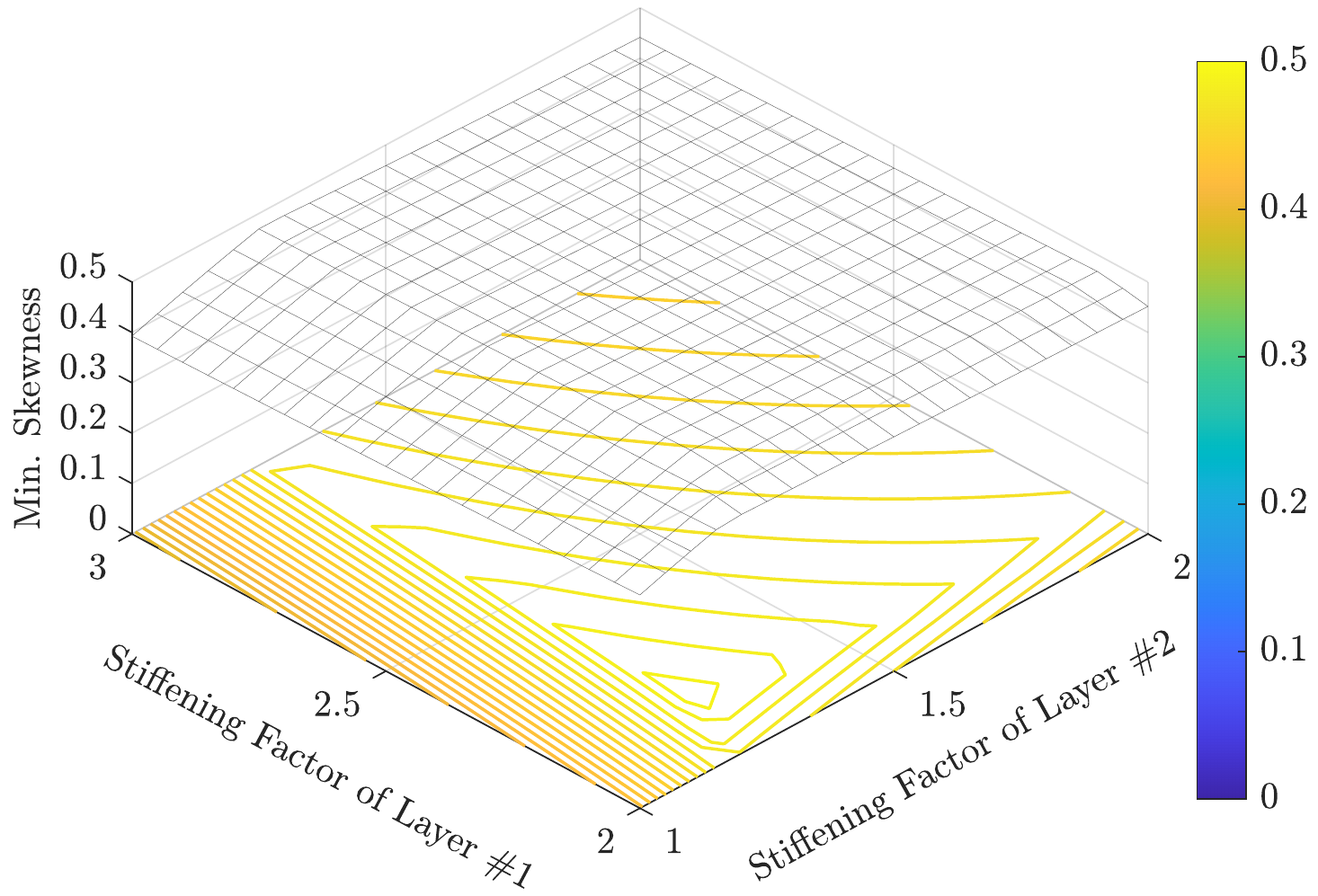}}
	\caption{\textcolor{black}{Effect of selective stiffening of two layers on the skewness metric of the beam in a channel test problem using the spring analogy model.}}
	\label{fig:SrngAnlg_sel_stiff_two_layers}
\end{figure*}

As for the \textbf{selective stiffening of two layers}, Fig. \ref{fig:SrngAnlg_sel_stiff_two_layers} presents the minimum skewness for the beam in a channel test problem using both diagonals simultaneously. First results are generated for stiffening factors ranging from 1 to 6 with a relatively large step of 0.1 as seen in Fig. \ref{figs:SrngAnlg_sel_stiff_two_layers_big_step}. The best combination of stiffening factors is in the vicinity of 2.2 and 1.3 for the first and second layers respectively, and produces a minimum skewness of 0.489. In Fig. \ref{figs:SrngAnlg_sel_stiff_two_layers_small_step}, a more detailed study is performed with a smaller step of 0.05. The minimum skewness attained is 0.493 at stiffening factors of 2.15 and 1.25 for the first and second layers respectively. Note that this improvement in minimum skewness is in comparison to a value of 0.327 for selective diagonals without stiffening (cf. Table \ref{tab:triangulation}) and a value of 0.432 for both diagonals with one layer stiffening (cf. Fig. \ref{figs:spring_sel_stiff_one_layer_skew}). The minimum skewness results of the remaining test problems are included in Table \ref{tab:summary_skew}.

This selective stiffening technique could be potentially extended to more layers at the fluid-structure interface. However, with more layers the search for the optimum values of the stiffening factors becomes even harder and the potential benefit becomes less significant.

This concludes our discussion on the spring analogy model, next we move on to discuss the second legacy mesh deformation technique; the continuous linear elasticity model.

\subsection{Continuous Linear Elasticity Mesh Deformation Model}
\label{ssec:linearelasticity}
\textcolor{black}{One of the obvious and relatively simple approaches to mesh deformation is modeling the mesh as a fictitious continuous elastic media, so it deforms according to the laws of linear elasticity.} \citet{Tezduyar1992computation} utilized a modified form of the linear homogeneous elasticity model to deform the fluid mesh. The modification includes dropping the Jacobian determinant in the calculation of the stiffness matrix, thus altering the transformation from the natural coordinates to the global coordinates. This results in introducing variable material properties where smaller elements have higher rigidity than larger ones, hence they are less prone to vanish during mesh deformation. This is potentially beneficial for boundary layers where smaller elements are usually used at the fluid-structure interface. In addition, they also employed selective treatment through altering the Lam\'e constants in order to stiffen elements more against shape changes rather than against volume changes. Later, \citet{Stein2003} introduced a new approach of mesh deformation techniques by varying the stiffening power to control the degree by which smaller elements are stiffened in comparison to larger ones. \citet{Stein2004} added more control to boundary layers by \say{glueing} the layers directly surrounding a moving, non-deformable solid objects (i.e., rigid body motion). Later, \citet{Jenkins2015, Jenkins2016} utilized this method to deform the fluid mesh in their TOFSI using a level set approach. Even though the before-mentioned mesh deformation techniques were applied mainly to unstructured, triangular meshing, \citet{Jenkins2015, Jenkins2016} utilized this technique on quad meshing.

\textcolor{black}{In the following discussion, we utilize the following mesh deformation problem statement:
\begin{equation}
	\begin{aligned}
		\mathbf{K}_{els} \mathbf{x} & = \mathbf{0} \quad \mathrm{on} \quad \mathrm{\Omega} \backslash \mathrm{\Gamma}_{\mathrm{FSI}}, \\
		\mathbf{x} & = \mathbf{u} \quad \mathrm{on} \quad \mathrm{\Gamma}_{\mathrm{FSI}}.
	\end{aligned}
\end{equation}
\noindent where $\mathbf{K}_{els}$ is the linear elastic stiffness matrix, $\mathbf{x}$ is the mesh nodal displacements vector, $\mathrm{\Omega}$ is the mesh computational domain, $\mathrm{\Gamma_{\mathrm{FSI}}}$ is the fluid-structure interface, and $\mathbf{u}$ is the structural nodal displacements vector.}

\subsubsection{Layered Selective Stiffening at the Fluid-Structure Interface}
\label{secss:selective_stiffening_elstcty}
\begin{figure*}
	\centering
	\subfloat{\includegraphics[height=0.02\linewidth]{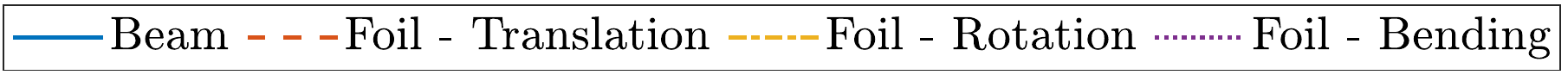}} \\
	\addtocounter{subfigure}{-1}
	\subfloat[\label{figs:LinEls_OneLayer_Skw}]{\includegraphics[width=0.33\linewidth]{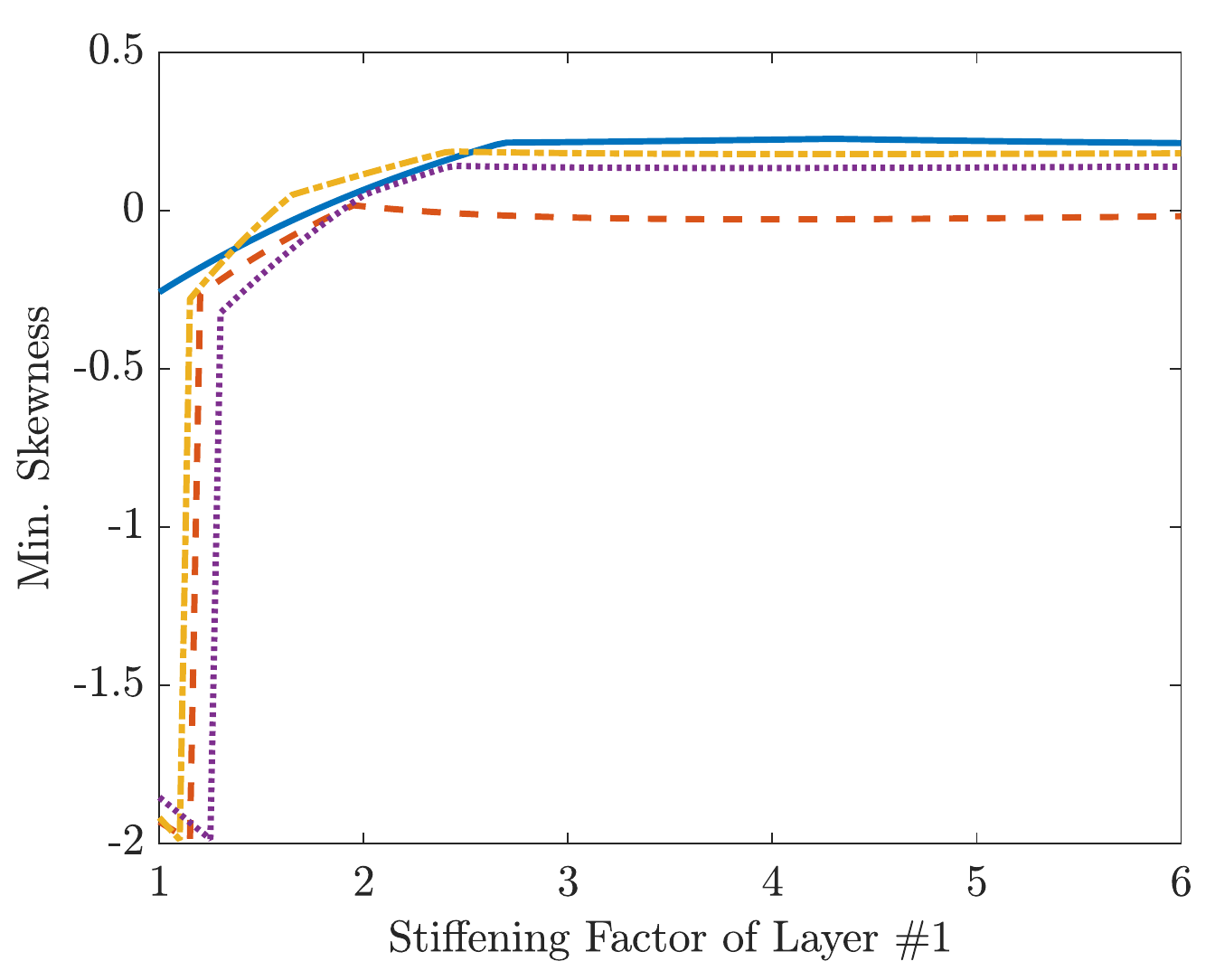}}
	\subfloat[]{\includegraphics[width=0.33\linewidth]{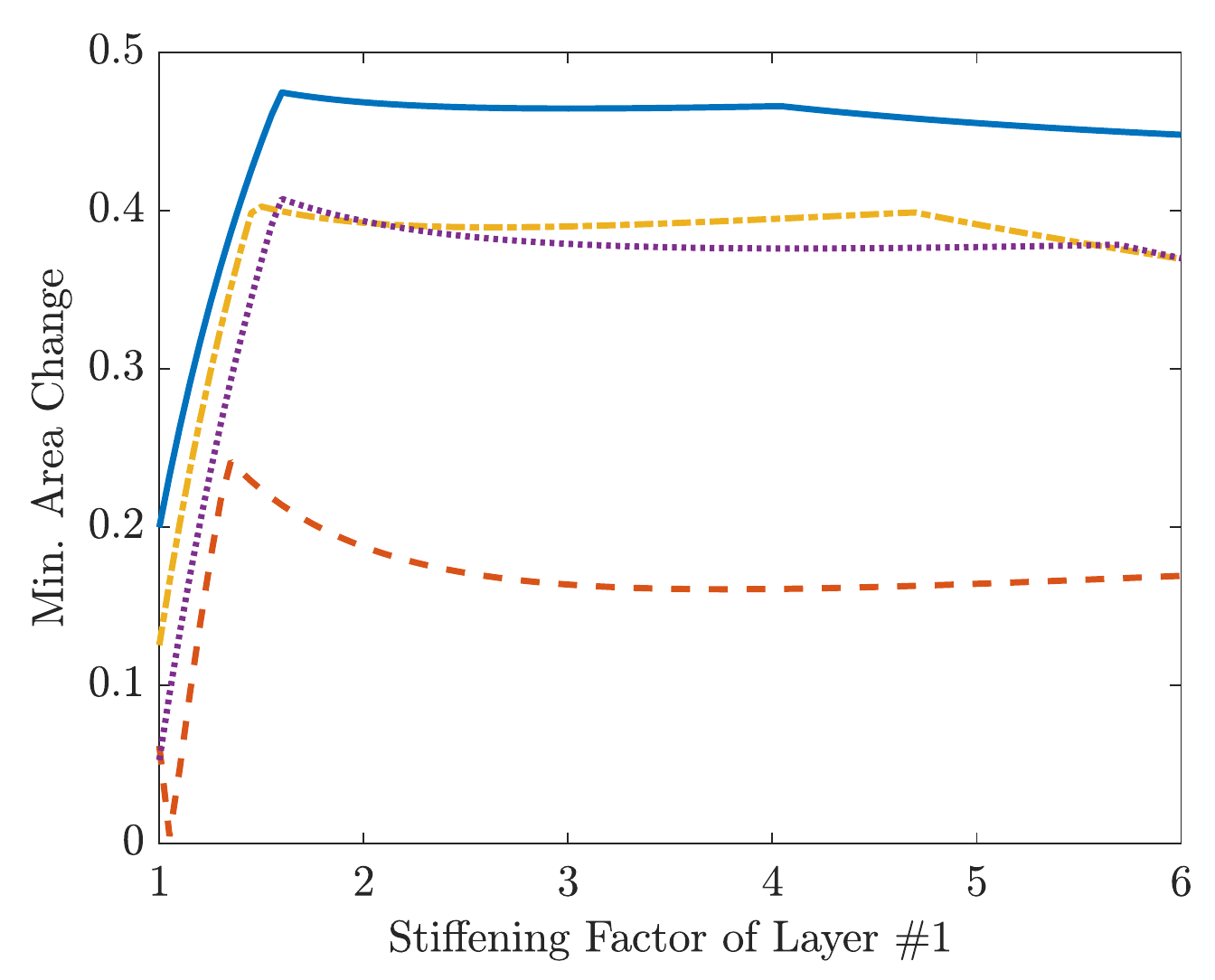}}
	\subfloat[]{\includegraphics[width=0.33\linewidth]{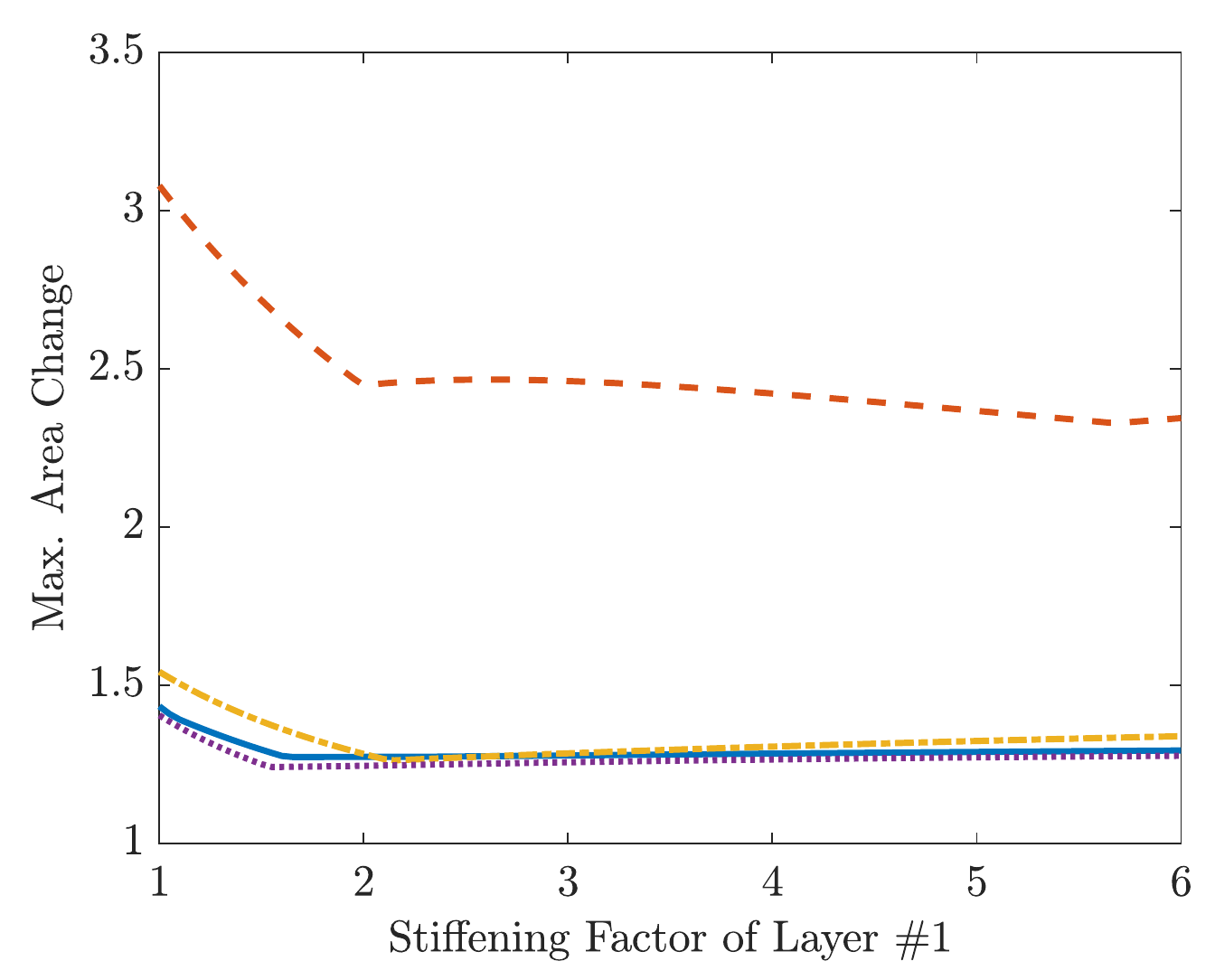}}
	\caption{Effect of selective stiffening of one layer on the mesh quality of the test problems using the linear elasticity model.}
	\label{fig:LinEls_sel_stiff_one_layer}
\end{figure*}

In this work, we implement a \textit{layered selective stiffening} technique similar to \ref{secss:selective_stiffening_sprng}. This is performed simply by multiplying the local stiffness matrix by a stiffening factor, which is equivalent to modifying the elastic modulus. Figure \ref{fig:LinEls_sel_stiff_one_layer} shows the mesh quality metrics of the test problems using the linear elasticity model with \textbf{one layer selective stiffening}. It shows a range of stiffening factors from 1 to 6 with a step of 0.05. It can clearly be seen that the linear elasticity model performs much worse than the spring analogy model even with one layer selective stiffening. The foil in a channel under translation has inverted elements for almost all values of the stiffening factor. The other test problems only start to not have inverted elements after a certain value of the stiffening factor (1.5 to 1.8). It seems that this model would greatly benefit from selective stiffening of more layers.

\begin{figure*}
	\centering
	\subfloat[\label{figs:LinEls_TwoLayers_Big}]{\includegraphics[width=0.495\linewidth]{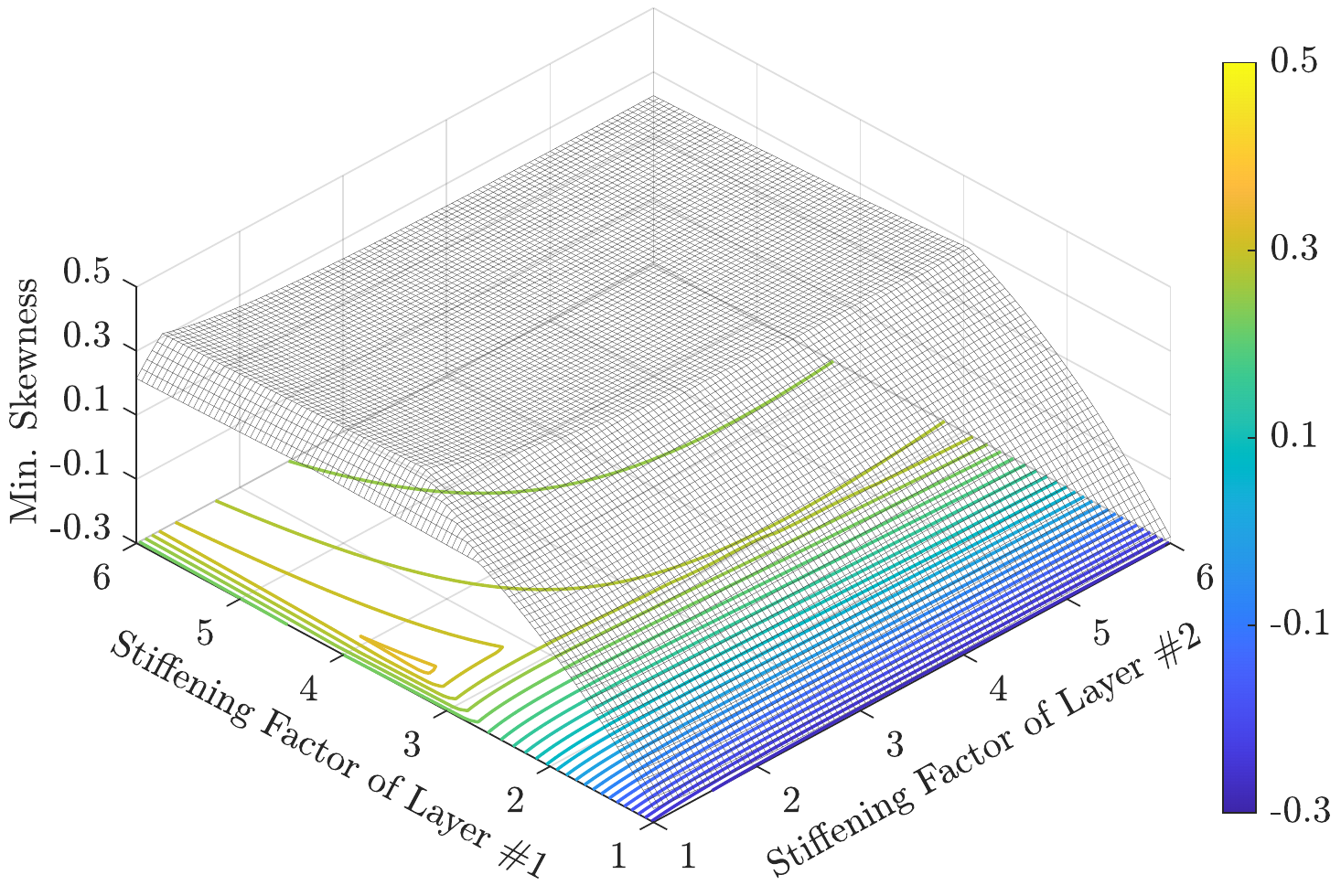}}
	\subfloat[\label{figs:LinEls_TwoLayers_Small}]{\includegraphics[width=0.495\linewidth]{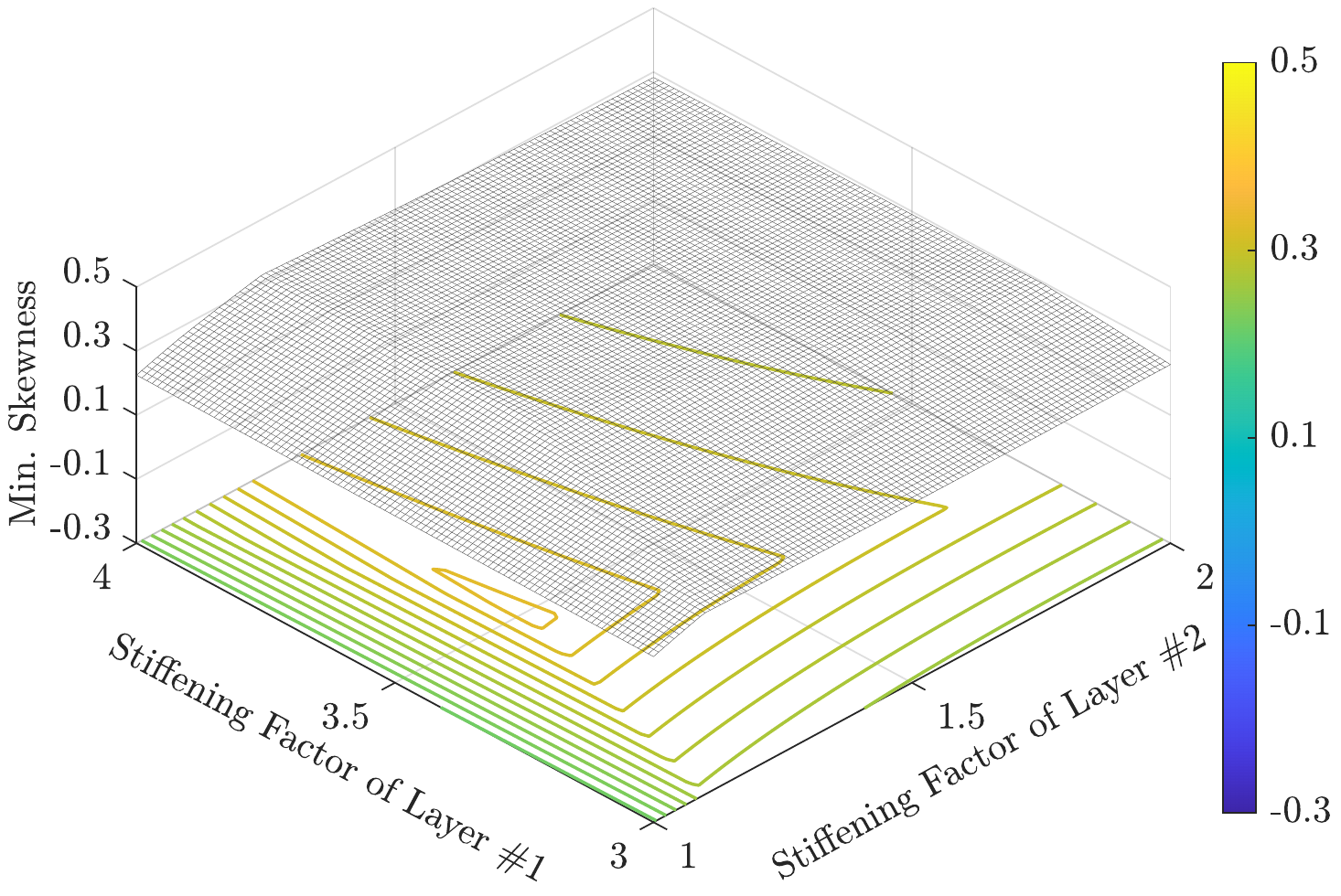}}
	\caption{\textcolor{black}{Effect of selective stiffening of two layers on the skewness metric of the beam in a channel test problem using the linear elasticity model.}}
\end{figure*}

Figure \ref{figs:LinEls_TwoLayers_Big} shows the \textbf{selective stiffening of two layers} of the beam in a channel problem. The range for both factors is 1 to 6 with a step of 0.05. Inverted elements with negative skewness are predominantly located in the area with a low stiffening factor for the first layer. The best minimum skewness attained is 0.332 at stiffening factors of 3.5 and 1.25 for the first and second layers respectively. A more detailed study with a smaller step of 0.01 is shown in Fig. \ref{figs:LinEls_TwoLayers_Small}. A slightly higher minimum skewness is located at 3.48 and 1.25 first and second stiffness factors respectively. In hindsight, a step of 0.05 for the two layers selective stiffening is more than enough for this study. Note that this result of the beam in a channel test problem is in comparison to a minimum skewness of -0.257 without stiffening and a minimum skewness of 0.227 with one layer selective stiffening factor of 4.3 (cf. Fig. \ref{figs:LinEls_OneLayer_Skw} for both values). Table \ref{tab:summary_skew} includes the minimum skewness results for all test problems upon using two layers selective stiffening. Three of the test problems achieve their best minimum skewness around the same values of the stiffening factors. The foil in a channel under translation test problem - which behaved badly with one layer selective stiffening - seems to behave well with two layers stiffening, although at rather high factors compared to the other problems.

As for the \textbf{selective stiffening of three layers}, the minimal computational costs means it's affordable to perform such calculations. A study is performed on all test problems with a stiffness factor range of 1 to 6 with a step of 0.1. Table \ref{tab:summary_skew} includes the minimum skewness results obtained and the values of the corresponding stiffness factors. The foil in a channel problem under translation in particular had its best minimum skewness at the limit of the stiffness factors range 1 to 6, so we had to widen the range. While the beam in a channel test problem shows a minor increase ($\approx$13\%) in the best minimum skewness with three layers selective stiffening, the foil in a channel test problems show a much better improvement ($\approx$25\% to 30\%).

A final remark about the continuous linear elasticity model is in order. \citet{Shamanskiy2021} mentioned that a better result could be obtained from this model if it is used in iterations instead of in a single step.

\section{Hyperelastic Model}
\textcolor{black}{A more advanced approach compared to the linear elasticity model is to represent the mesh as a hyperelastic media since hyperelasticity models are better suited to model shape deformations such as skewing. The same conditions describing the linear elasticity model problem apply for the hyperelastic model except that the material laws are different.} \citet{Dettmer2006} utilized a pseudo-elastic model based on a simple hyperelastic Neo-Hookean material to represent the fluid  mesh. Then, they used a simplified strategy to calculate the nodal positions so as to optimize the mesh quality based on \citet{Braess2000}. \citet{Takizawa2020} used fiber-reinforced hyperelasticity with optimized zero-stress states to develop a low-distortion mesh moving method. More recently, \citet{Shamanskiy2021} introduced a mesh deformation technique based on a logarithmic variation of the neo-Hookean material model.

In this work, we suggest using the Yeoh hyperelastic material model as a mesh deformation technique. Yeoh's model was originally introduced in \citep{Yeoh1993} as an improvement over the classic hyperelasticity models such as the neo-Hookean and Mooney-Rivlin for modeling elastomers. In its nearly incompressible form, the Yeoh strain energy function is described as follows:
\begin{equation}
	\begin{aligned}
			W = & A_{10} (J_1 - 3) + A_{20} (J_1 - 3)^2 + A_{30} (J_1 - 3)^3 \\
			& + \kappa (J_3 - 1)^2. 
	\end{aligned}
\end{equation}

\noindent where $A_{10}$, $A_{20}$, and $A_{30}$ are material parameters, $J_1$ and $J_3$ are the first and third reduced strain invariants, and $\kappa$ is the bulk modulus upon which depends the degree of incompressibility imposed. \citet{Yeoh1993} mentioned that his proposed model is conceptually a material model with a shear modulus that varies with deformation. And since the main stress mode that affects element skewness is shear, it makes sense that this model behaves well in this aspect. Unlike a linear elastic model where the value of the elastic modulus is completely arbitrary\footnote{Probably the only consideration in this case is to ensure that the mathematical system is well-scaled so as to avoid any numerical issues due to differing orders of magnitude.}, in hyperelasticity the values of the material constants do have an effect on the quality of the output mesh and to some degree the convergence of the problem.

In this work, we employ the total Langrangian formulation along with the second Piola-Kirchhoff stress and the Lagrangian strain definitions. \textcolor{black}{The problem statement for this mesh deformation model is as follows:
\begin{equation}
	\begin{aligned}
		{}^t\mathbf{K} \; \Delta \mathbf{x} & = {}^{ex}\mathbf{F} - {}^{in}\mathbf{F} \quad \mathrm{on} \quad \mathrm{\Omega}, \\
		\mathbf{x} & = \mathbf{u} \quad \quad \quad \quad \, \, \, \, \, \mathrm{on} \quad \mathrm{\Gamma_{\mathrm{FSI}}}.
	\end{aligned}
\end{equation}
\noindent where ${}^t\mathbf{K}$ is the tangent stiffness matrix, $\mathbf{x}$ is the incremental nodal displacements vector, ${}^{ex}\mathbf{F}$ is the external nodal forces vector, ${}^{in}\mathbf{F}$ is the internal nodal forces vector, $\mathrm{\Omega}$ is the mesh computational domain, $\mathbf{u}$ is the structural nodal displacements vector, and $\mathrm{\Gamma_{\mathrm{FSI}}}$ is the fluid-structure interface. This problem is to be solved iteratively (e.g., Newton-Raphson), where the prescribed displacements at the fluid-structure interface are increased incrementally until convergence.} For details on the finite element implementation, see \citep[p.~492]{Bathe2014finite} or \citep[p.~200]{kim2014introduction}. In the following, we study in detail some aspects of the hyperelastic model as a mesh deformation technique.

\subsection{Parametric Study on the Yeoh Model}
From numerical experiments, the material parameters $A_{10}$ and $A_{30}$ don't have a considerable effect on the mesh quality, hence they are set to 1 and 0 respectively. The parametric study in this subsection is focused on the material parameter $A_{20}$ and the bulk modulus $\kappa$.

\begin{figure*}
	\centering
	\subfloat{\includegraphics[height=0.02\linewidth]{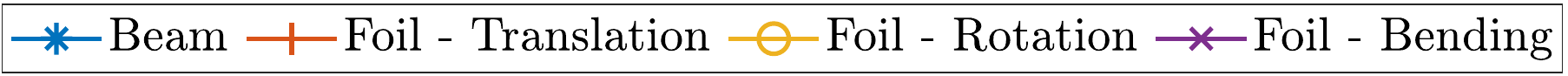}} \\
	\addtocounter{subfigure}{-1}
	\subfloat[]{\includegraphics[width=0.245\linewidth]{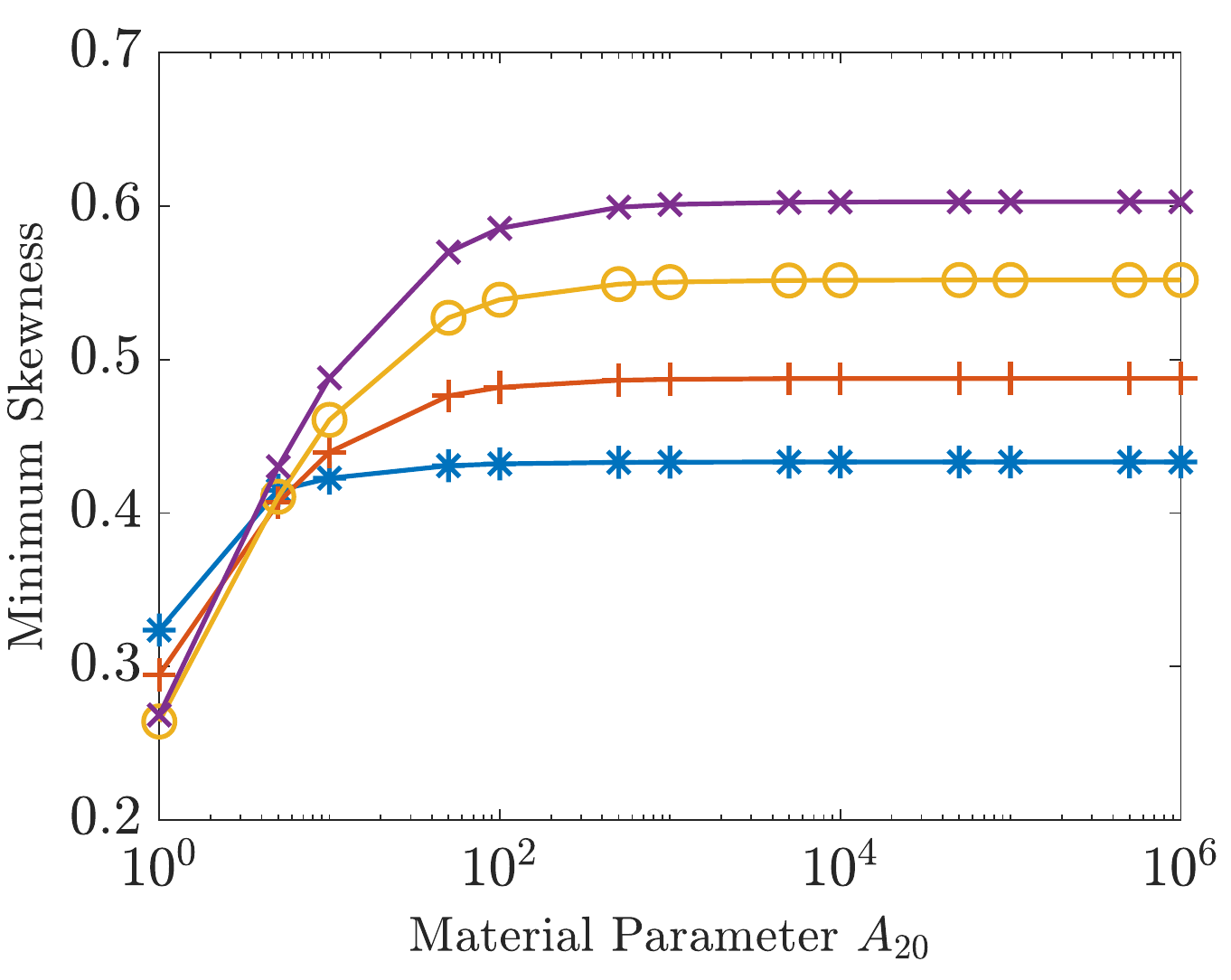}}
	\subfloat[]{\includegraphics[width=0.245\linewidth]{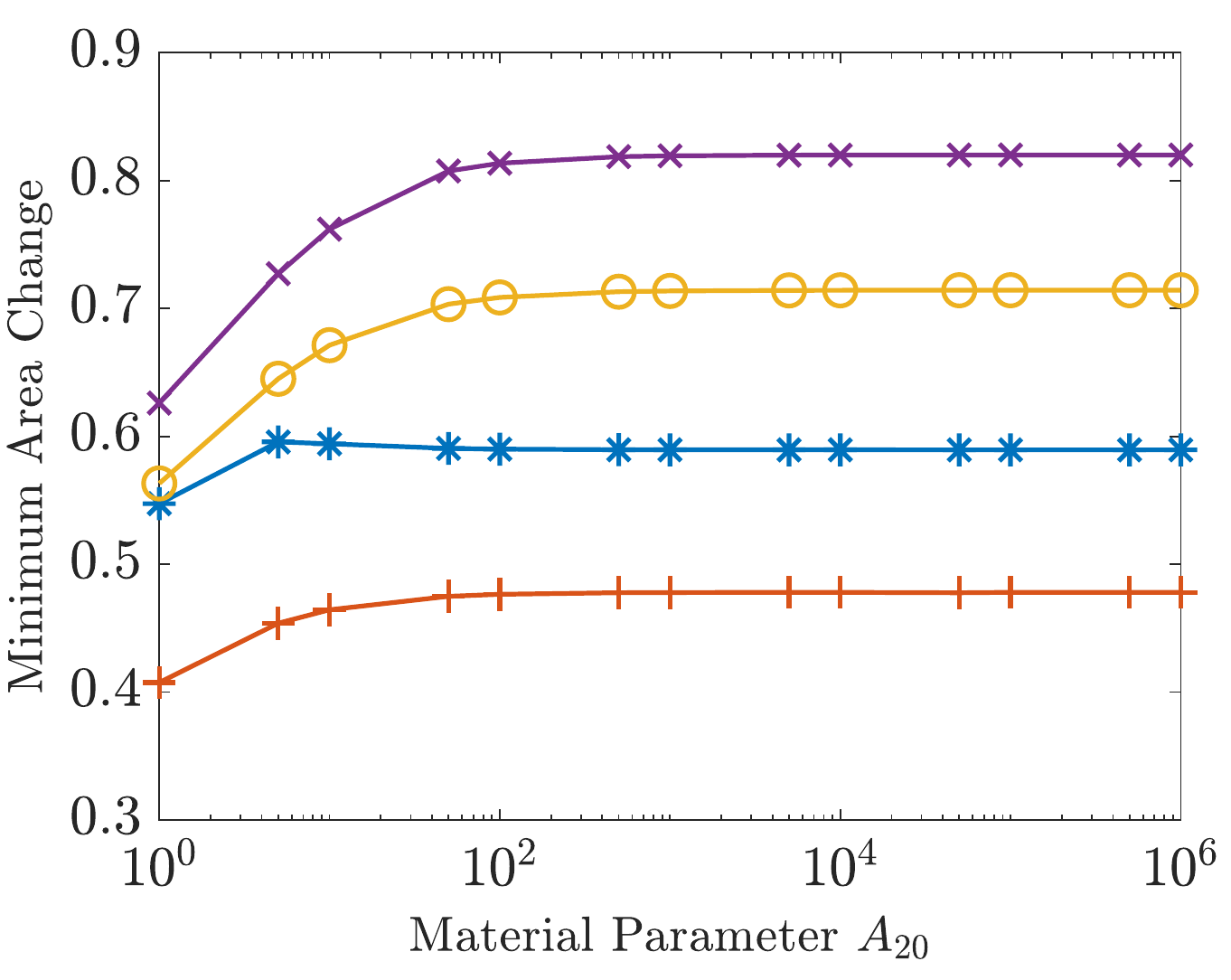}} 
	\subfloat[]{\includegraphics[width=0.245\linewidth]{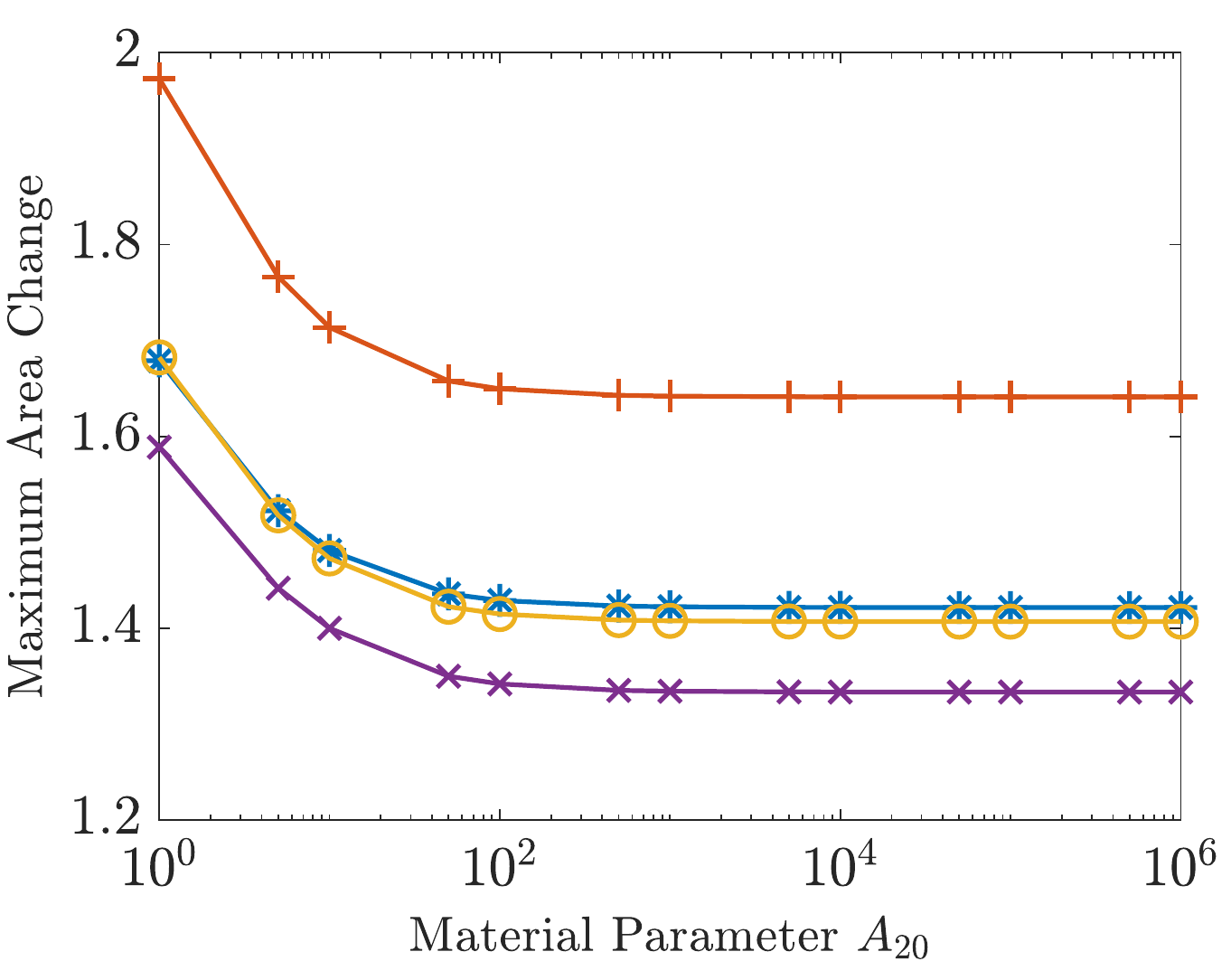}}
	\subfloat[\label{figs:effect_of_a20_time}]{\includegraphics[width=0.245\linewidth]{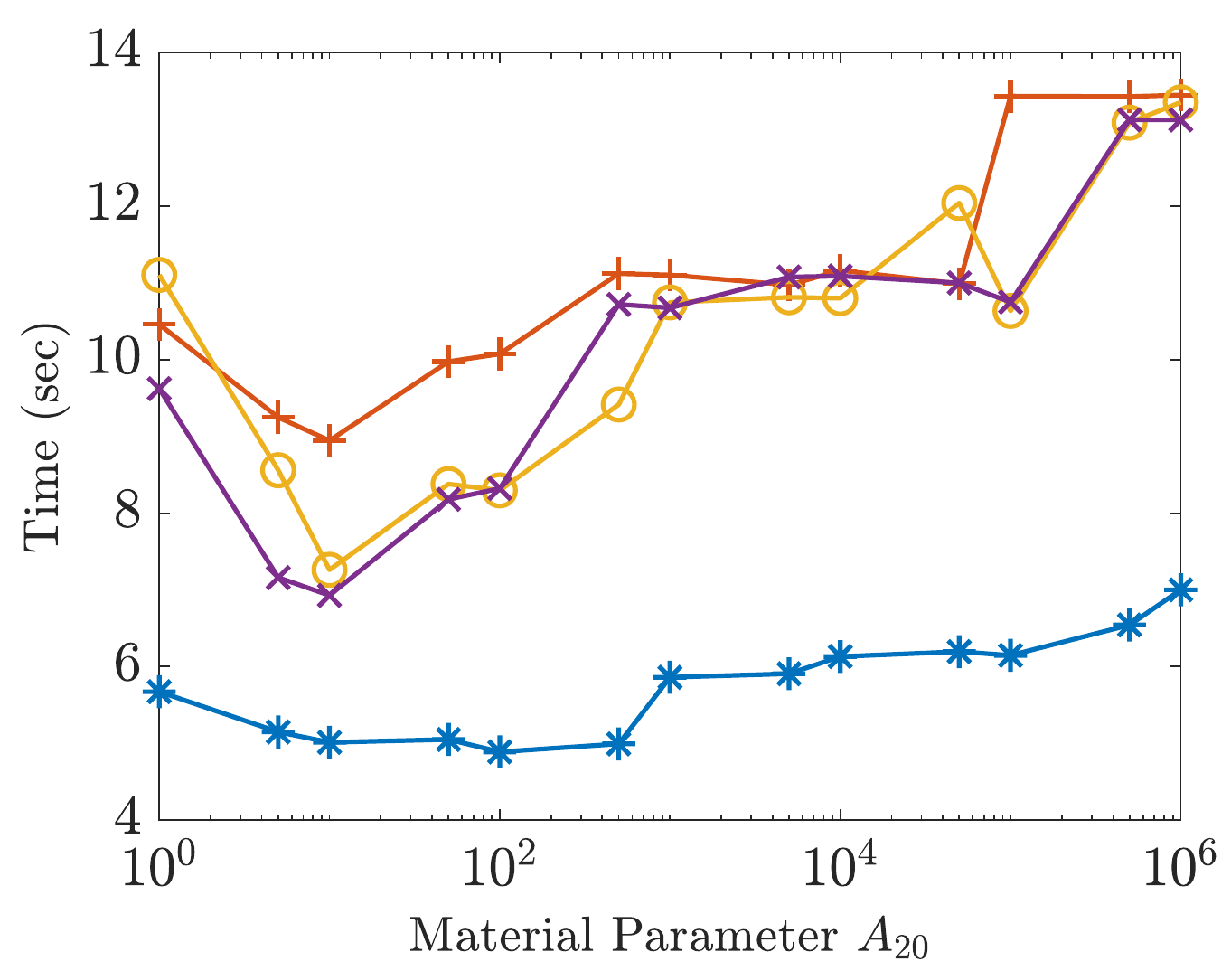}}
	\caption{Effect of the material parameter $A_{20}$ on the mesh quality of the test problems using the Yeoh hyperelastic model.}
	\label{fig:effect_of_a20}
\end{figure*}

As for the effect of the material parameter $A_{20}$ on the mesh quality, we solved all test problems using different values of $A_{20}$ and the results are shown in Fig. \ref{fig:effect_of_a20}. The bulk modulus $\kappa$ in this particular study is set to 1. It seems that the best mesh metrics - for $\kappa = 1$ - are obtained around a value of $A_{20} = 10^3$, and they plateau for higher values of $A_{20}$. For the solution time (Fig. \ref{figs:effect_of_a20_time}), a value of $A_{20} = 10^3$ doesn't give the best computational time but it's not much higher than the near optimum computational time obtained at $A_{20} = 10^1$.

\begin{figure*}
	\centering
	\subfloat{\includegraphics[height=0.02\linewidth]{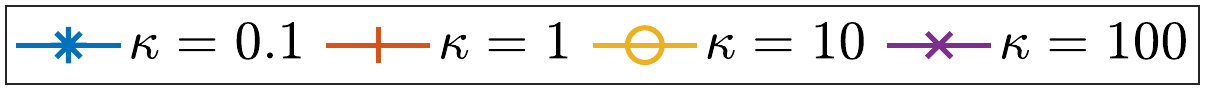}} \\
	\addtocounter{subfigure}{-1}
	\subfloat[\label{figs:effect_of_a&k_skw}]{\includegraphics[width=0.245\linewidth]{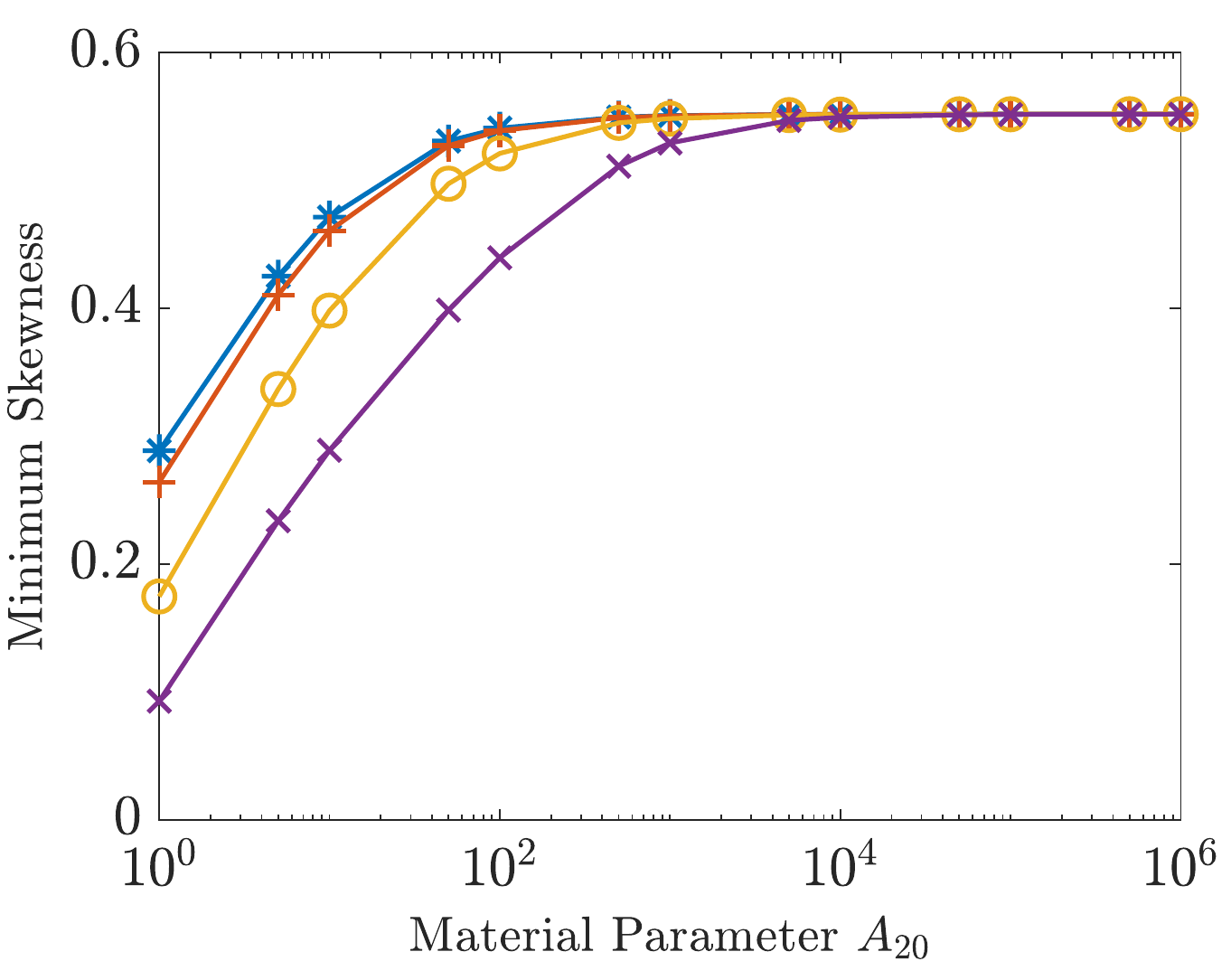}}
	\subfloat[\label{figs:effect_of_a&k_mina}]{\includegraphics[width=0.245\linewidth]{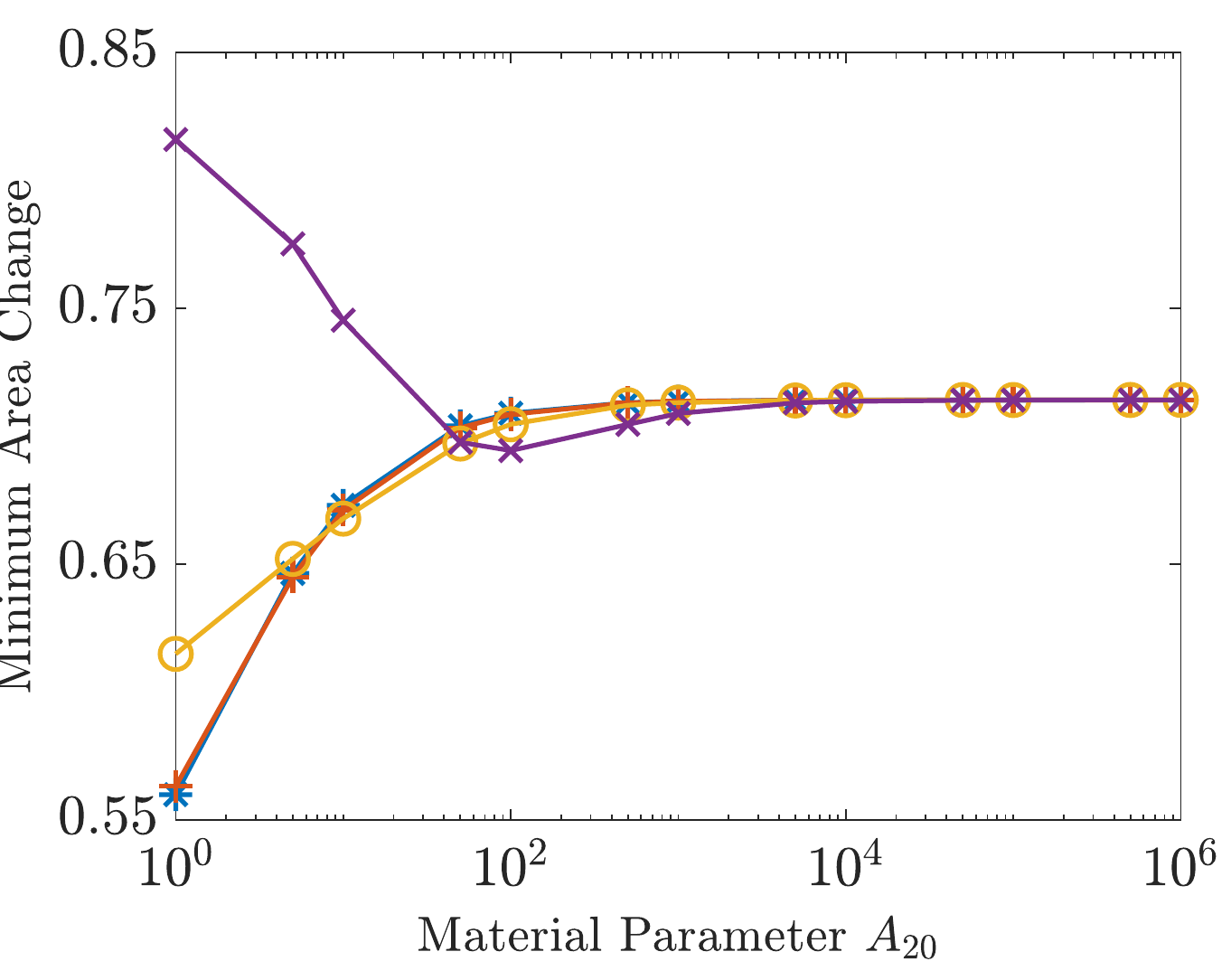}} 
	\subfloat[\label{figs:effect_of_a&k_maxa}]{\includegraphics[width=0.245\linewidth]{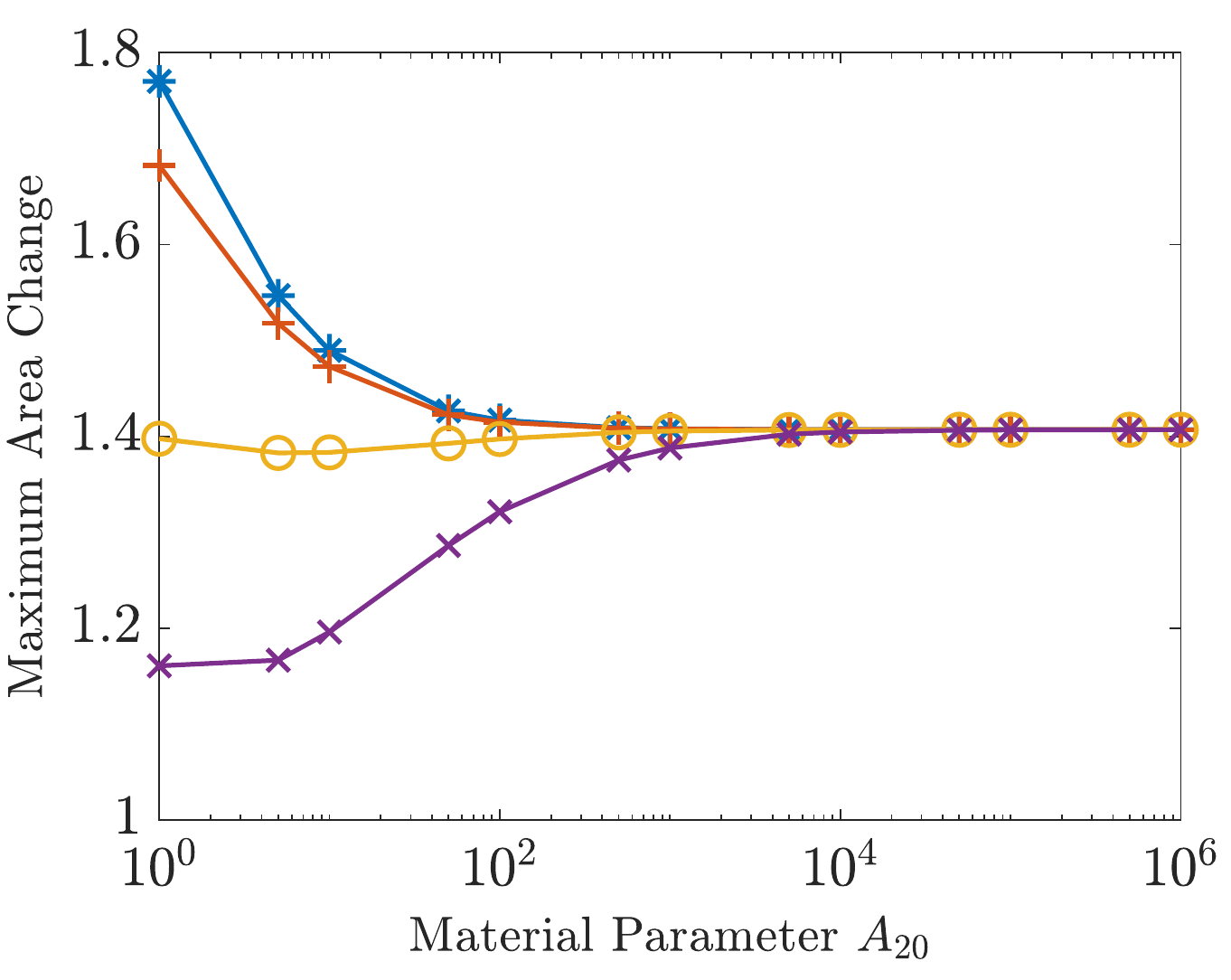}}
	\subfloat[\label{figs:effect_of_a&k_tm}]{\includegraphics[width=0.245\linewidth]{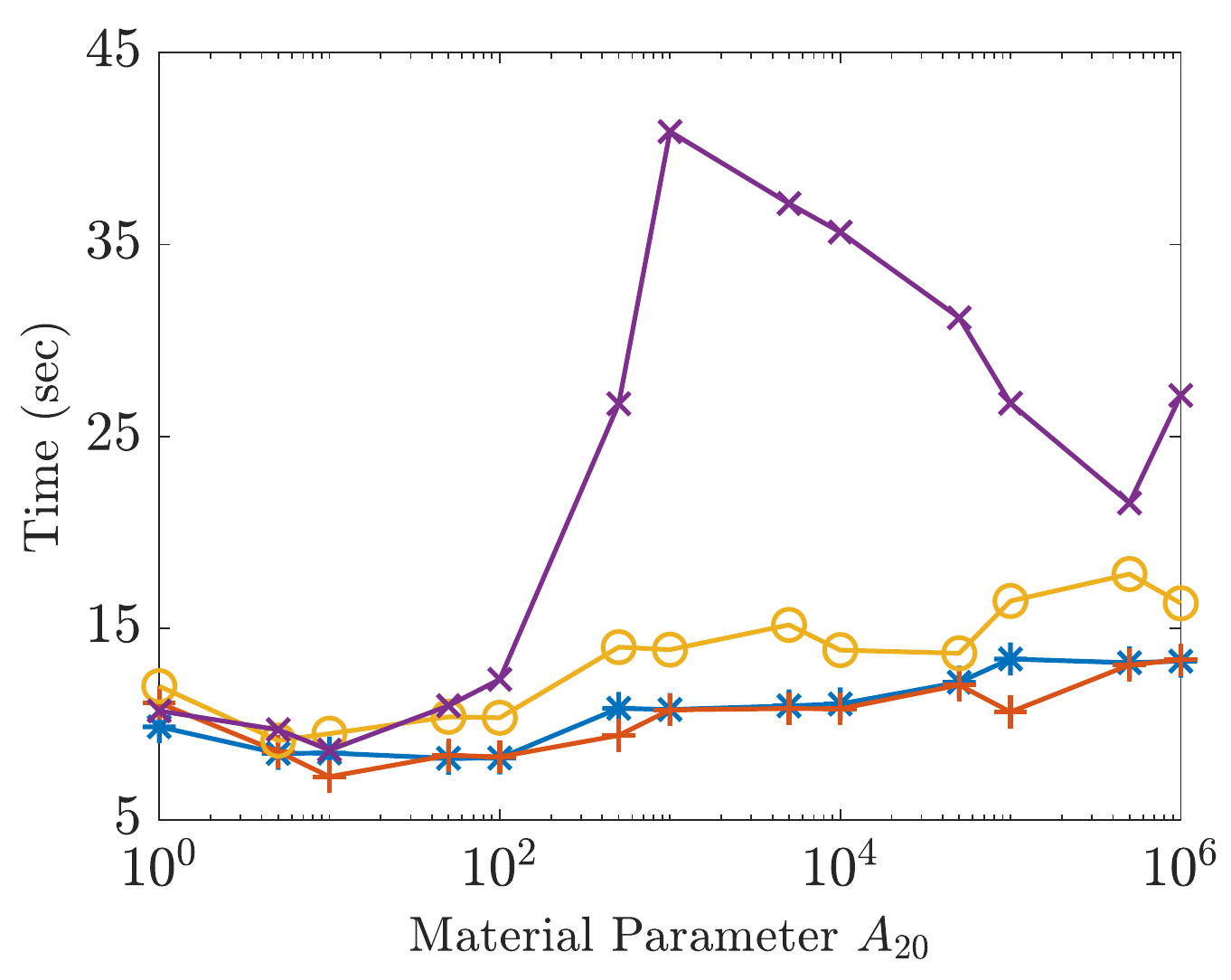}}
	\caption{Interplay of the material parameter $A_{20}$ and the bulk modulus $\kappa$ and their effect on the mesh quality of the foil in a channel under rotation test problem using the Yeoh hyperelastic model.}
	\label{fig:effect_of_a&k}
\end{figure*}

To investigate in more detail the interplay of the material parameter $A_{20}$ and the bulk modulus $\kappa$, we solved the foil in a channel under rotation test problem for a range of $A_{20}$ and different values of $\kappa$ (Fig. \ref{fig:effect_of_a&k}). As expected, increasing the bulk modulus $\kappa$ has a positive effect on the volume change mesh metrics (Figs. \ref{figs:effect_of_a&k_mina} \& \ref{figs:effect_of_a&k_maxa}), however this effect diminishes for values of the material parameter $A_{20}$ higher than $10^2$. On the other hand, increasing $\kappa$ has a negative effect on the minimum skewness (Fig. \ref{figs:effect_of_a&k_skw}) and the computational time (Fig. \ref{figs:effect_of_a&k_tm}). In addition, decreasing $\kappa$ further than 1 doesn't have that big an effect as evident in the close overlap of the data of $\kappa=0.1$ and $\kappa=1$ in Fig. \ref{fig:effect_of_a&k}.

\subsection{Layered Selective Stiffening at the Fluid-Structure Interface}
\begin{figure*}
	\centering
	\subfloat{\includegraphics[height=0.02\linewidth]{Fig_horizontal_legend_nomarkers}} \\
	\addtocounter{subfigure}{-1}
	\subfloat[\label{figs:hyper_twolayers_skw}]{\includegraphics[width=0.33\linewidth]{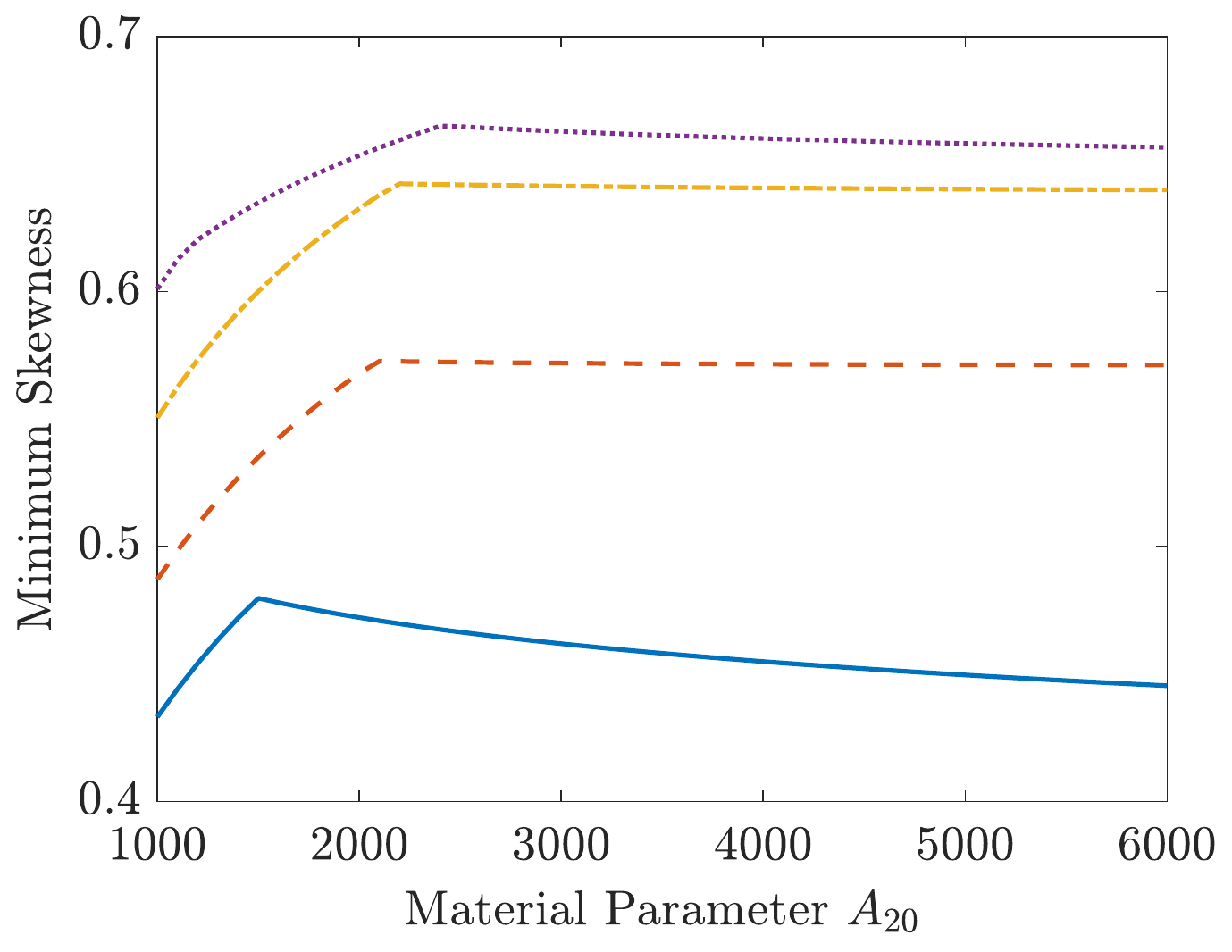}}
	\subfloat[\label{figs:hyper_twolayers_mina}]{\includegraphics[width=0.33\linewidth]{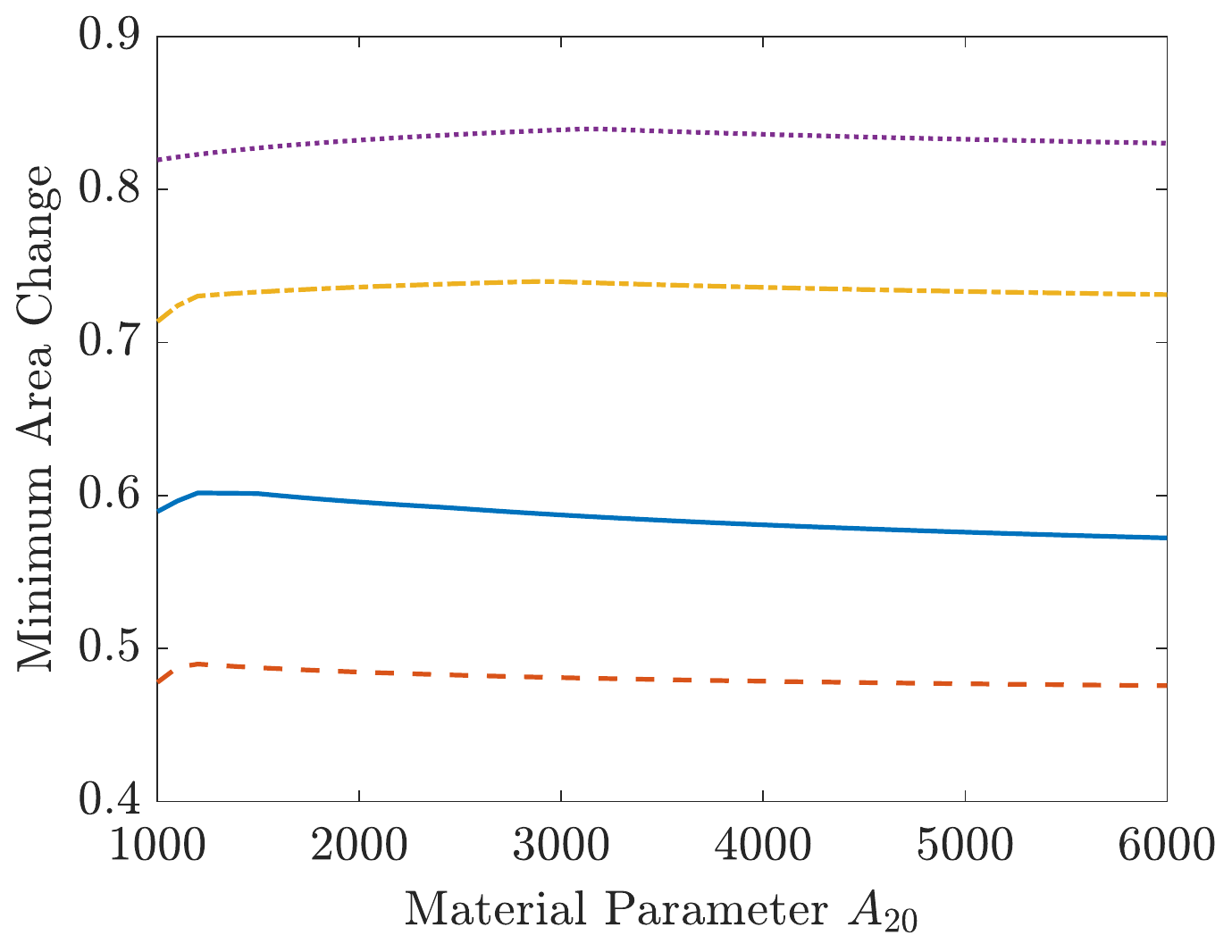}}
	\subfloat[\label{figs:hyper_twolayers_maxa}]{\includegraphics[width=0.33\linewidth]{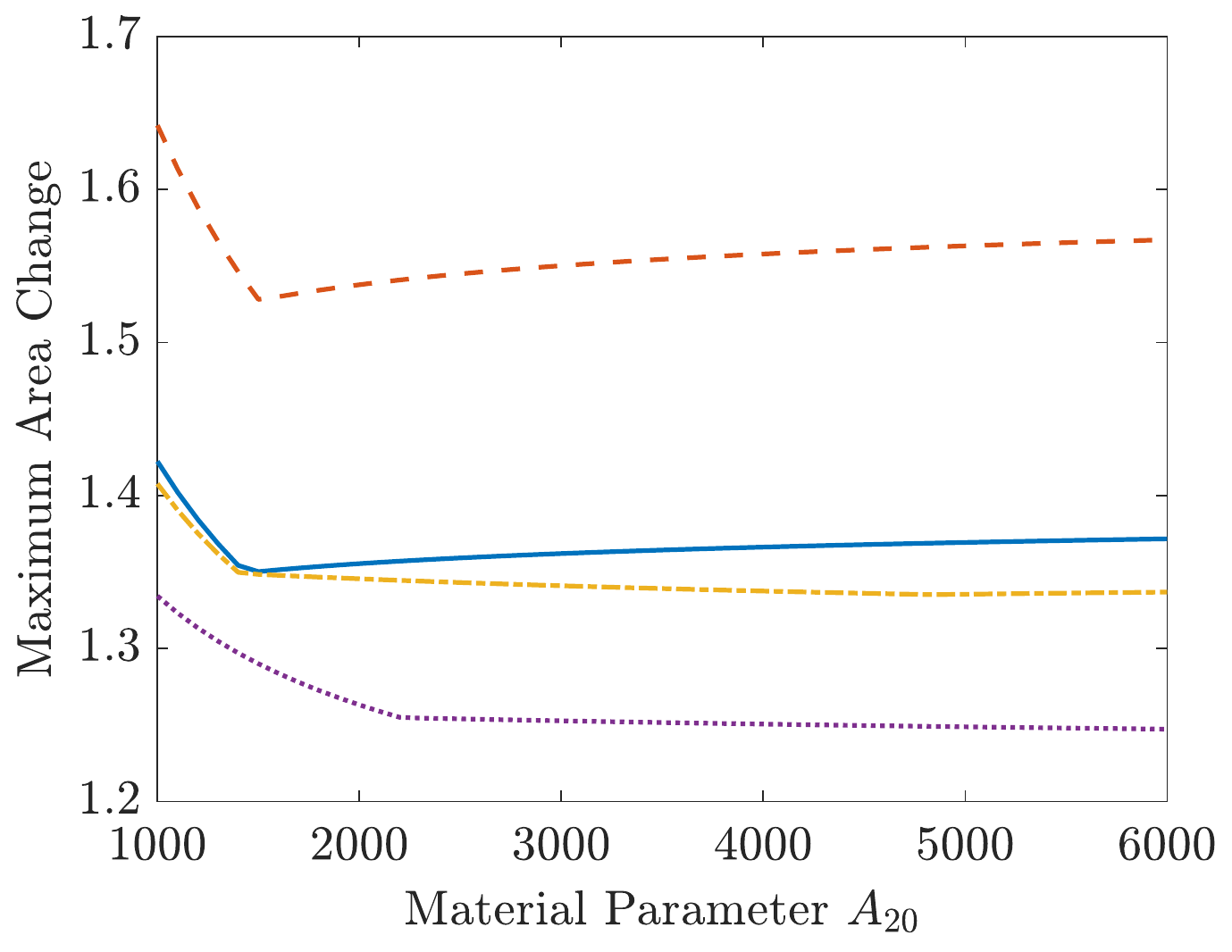}}
	\caption{Effect of selective stiffening of one layer on the mesh quality of the test problems using the hyperelasticity model.}
	\label{fig:Hyper_sel_stiff_one_layer}
\end{figure*}

Similar to subsections \ref{secss:selective_stiffening_sprng} and \ref{secss:selective_stiffening_elstcty}, we implement \textit{layered selective stiffening} with the hyperelastic model, where the material parameter $A_{20}$ - which is set to $10^3$ - is multiplied by a stiffening factor. Figure \ref{fig:Hyper_sel_stiff_one_layer} shows the effect of \textbf{one layer selective stiffening} on the mesh metrics of the test problem. The one layer selective stiffening increases the minimum skewness (Fig. \ref{figs:hyper_twolayers_skw}) even higher, which further proves the feasibility of the hyperelastic model as a mesh deformation technique. The minimum skewness in the foil in a channel test problems seems to plateau after reaching its highest value, but this is not the case with the beam in a channel test problem. Increasing the stiffening factor doesn't seem to have an effect on the minimum area change (Fig. \ref{figs:hyper_twolayers_mina}), but it does have a positive effect on the maximum area change (Fig. \ref{figs:hyper_twolayers_maxa}).

\begin{figure*}
	\centering
	\subfloat[]{\includegraphics[width=0.495\linewidth]{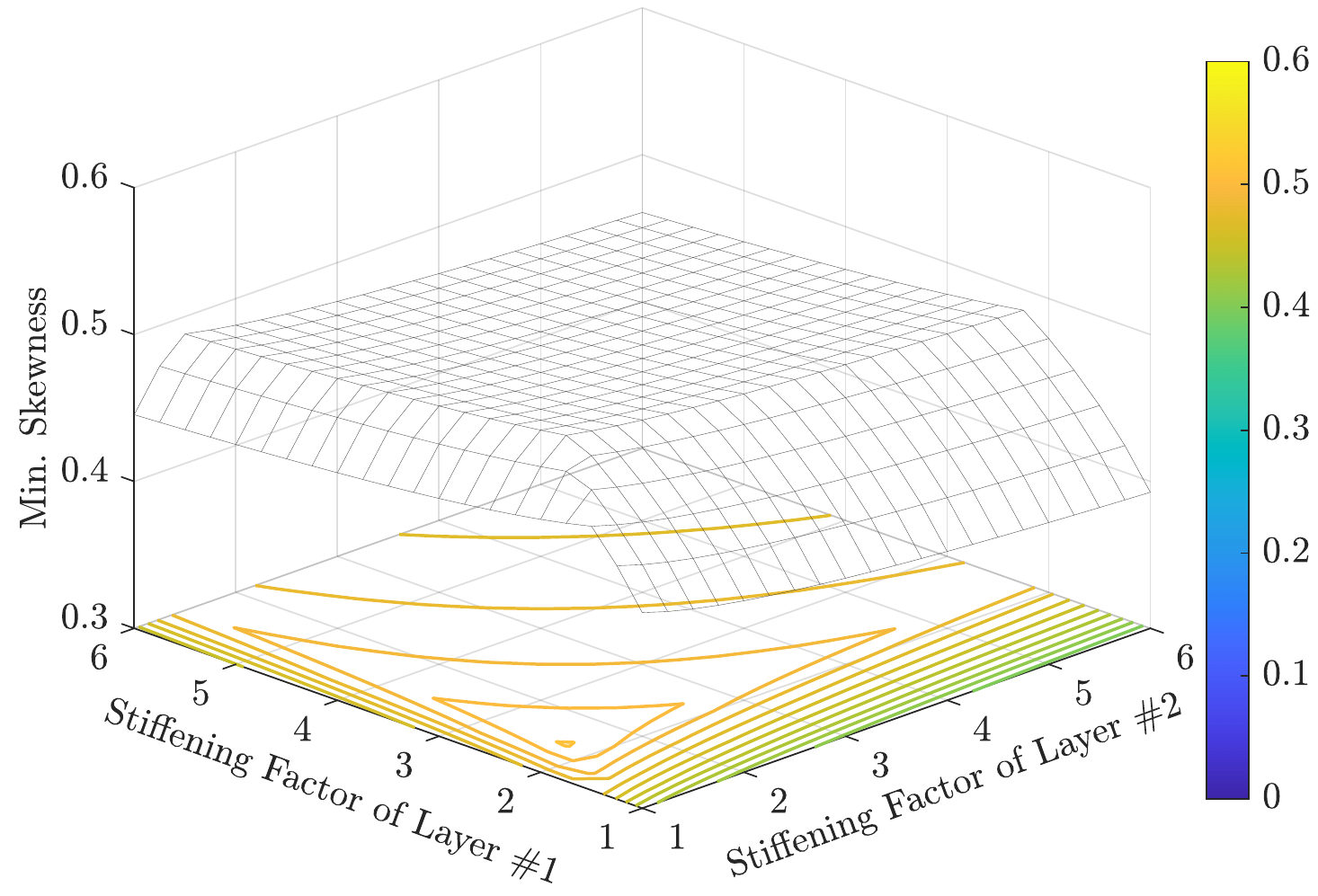}}
	\subfloat[]{\includegraphics[width=0.495\linewidth]{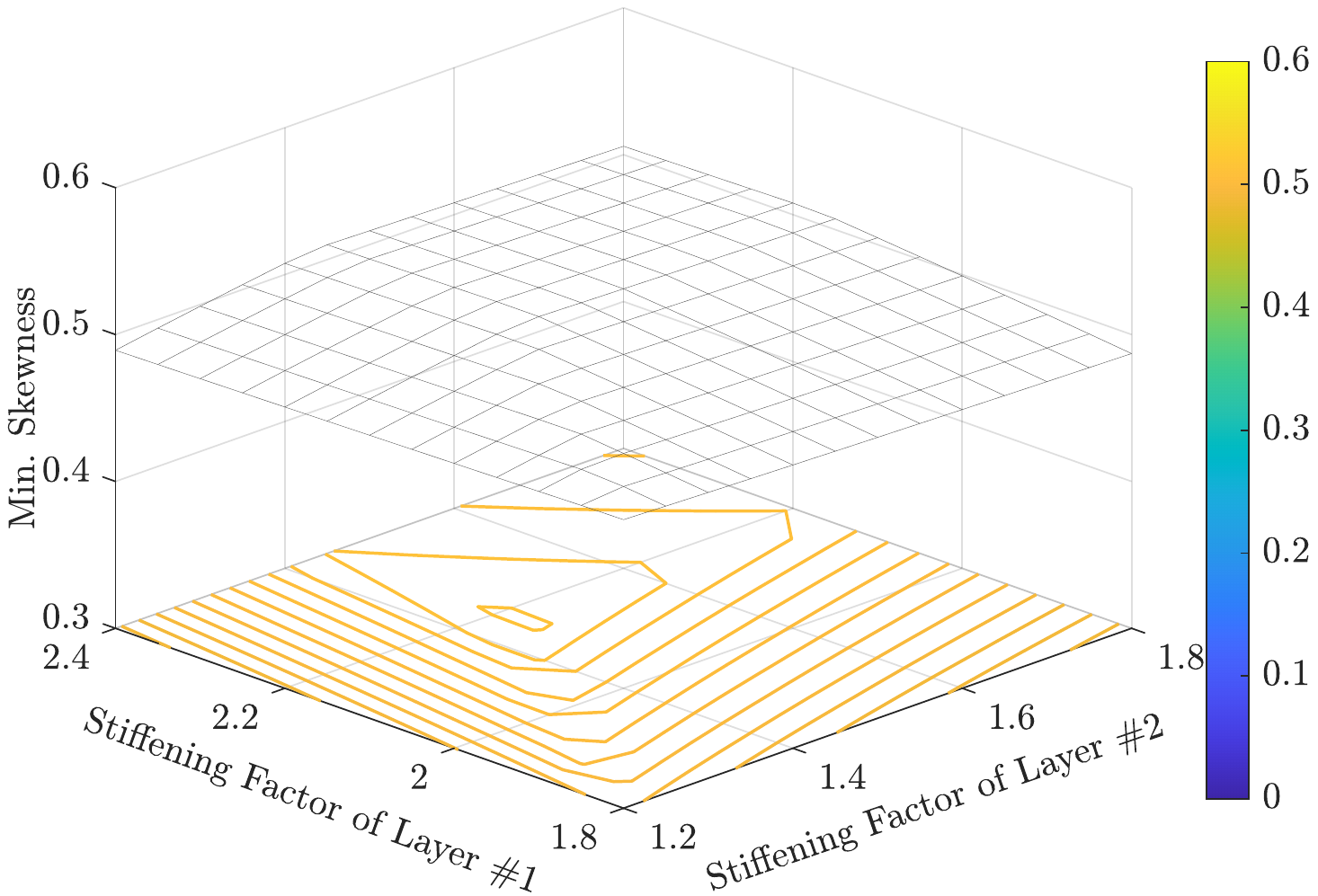}}
	\caption{\textcolor{black}{Effect of selective stiffening of two layers on the minimum skewness of the beam in a channel test problem using the hyperelasticity model.}}
	\label{fig:Beam_Hyper_sel_stiff}
\end{figure*}

Next, we utilized \textbf{two layers selective stiffening} on all the test problems using the hyperelastic model. We used a range of 1 to 6 for the stiffening factors with a relatively large step of 0.25 due to the high computational cost of solving the hyperelastic model. Once a good combination of stiffening factors is found, we tightened the search range and used a smaller step of 0.05. Even better results are obtained from using two layers selective stiffening. Figure \ref{fig:Beam_Hyper_sel_stiff} shows minimum skewness results of the beam in a channel test problem and the results of the remaining test problems are included in Table \ref{tab:summary_skew}.

In the next subsection, we discuss the sensitivity analysis of the hyperelastic mesh deformation model for implementation in topology optimization of fluid-structure interactions.

\begin{table*}[t!]
	\centering
	\small
	\caption{Summary of all minimum skewness results.}
	\label{tab:summary_skew}
	\begin{tabular}{@{}lllllllllll@{}} \toprule
		\multirow{5}{*}{Test Problem} & \multicolumn{3}{c}{Spring Analogy} & \multicolumn{4}{c}{Linear Elasticity} & \multicolumn{3}{c}{Hyperelastic} \\
		\cmidrule{2-11} 
		& \multicolumn{10}{c}{Minimum Skewness} \\ 
		& \multicolumn{10}{c}{\textcolor{red}{Stiffening Factor of Layer 1}} \\ 
		& \multicolumn{10}{c}{\textcolor{blue}{Stiffening Factor of Layer 2}} \\ 
		& \multicolumn{10}{c}{\textcolor{ForestGreen}{Stiffening Factor of Layer 3}} \\ \midrule
		Beam in a Channel & 0.327 & 0.432 & 0.493 & -0.257 & 0.227 & 0.332 & 0.377 & 0.433 & 0.480 & 0.512 \\
		& & \textcolor{red}{1.80} & \textcolor{red}{2.15} & & \textcolor{red}{4.30} & \textcolor{red}{3.50} & \textcolor{red}{4.00} & & \textcolor{red}{1.50} & \textcolor{red}{2.20} \\
		& & & \textcolor{blue}{1.25} & & & \textcolor{blue}{1.25} & \textcolor{blue}{1.60} & & & \textcolor{blue}{1.45} \\
		& & & & & & & \textcolor{ForestGreen}{1.20} & & \\
		Foil in a Channel - Translation & 0.338 & 0.420 & 0.506 & -1.930 & 0.017 & 0.369 & 0.470 & 0.488 & 0.573 & 0.613 \\
		& & \textcolor{red}{1.70} & \textcolor{red}{2.10} & & \textcolor{red}{1.95} & \textcolor{red}{5.40} & \textcolor{red}{7.10} & & \textcolor{red}{2.10} & \textcolor{red}{3.05} \\
		& & & \textcolor{blue}{1.30} & & & \textcolor{blue}{4.50} & \textcolor{blue}{3.20} & & & \textcolor{blue}{1.50}\\
		& & & & & & & \textcolor{ForestGreen}{1.80} & & \\
		Foil in a Channel - Rotation & 0.320 & 0.478 & 0.510 & -1.918 & 0.187 & 0.343 & 0.446 & 0.552 & 0.642 & 0.683 \\
		& & \textcolor{red}{4.10} & \textcolor{red}{2.50} & & \textcolor{red}{2.45} & \textcolor{red}{3.70} & \textcolor{red}{5.00} & & \textcolor{red}{2.20} & \textcolor{red}{3.30} \\
		& & & \textcolor{blue}{1.40} & & & \textcolor{blue}{1.70} & \textcolor{blue}{2.40} & & & \textcolor{blue}{1.55}\\
		& & & & & & & \textcolor{ForestGreen}{1.50} & & \\
		Foil in a Channel - Bending & 0.313 & 0.435 & 0.527 & -1.854 & 0.142 & 0.310 & 0.418 & 0.603 & 0.665 & 0.691 \\ 
		& & \textcolor{red}{3.00} & \textcolor{red}{2.50} & & \textcolor{red}{2.45} & \textcolor{red}{3.55} & \textcolor{red}{4.90} & & \textcolor{red}{2.40} & \textcolor{red}{3.75} \\ 
		& & & \textcolor{blue}{1.40} & & & \textcolor{blue}{1.60} & \textcolor{blue}{2.30} & & & \textcolor{blue}{1.45}\\
		& & & & & & & \textcolor{ForestGreen}{1.50} & & \\ \bottomrule
	\end{tabular} %
\end{table*}

\subsection{Sensitivity Analysis for Topology Optimization of Fluid-Structure Interactions}

In gradient-based topology optimization, it's a requirement to calculate the sensitivities of the fluid mesh deformation model in order to calculate the overall sensitivities of the optimization problem. Following the adjoint sensitivity of the three-field formulation presented in \citep{Maute2003}\footnote{\citet{Martins2005} discuss the adjoint-coupled sensitivity analysis in general.}, there are three sets of state variables in a TOFSI problem; namely \(\mathbf{u}\) the structural displacements, \(\mathbf{x}\) the fluid mesh nodal displacements, and \(\mathbf{w}\) the fluid velocities and pressures. For the adjoint sensitivity formulation, it's necessary to obtain the derivatives of the discrete equations governing the equilibrium of the solid, the fluid mesh deformation, and the fluid models - designated \(\mathbb{S}\), \(\mathbb{D}\), and \(\mathbb{F}\)  respectively - with respect to the three sets of state variables \(\mathbf{u}\), \(\mathbf{x}\), and \(\mathbf{w}\). We note that the governing equations of the fluid mesh motion \(\mathbb{D}\) are a function of only \(\mathbf{u}\) and \(\mathbf{x}\) but not \(\mathbf{w}\). The total derivative of a system variable $f$ (e.g., an objective function or a constraint) with respect to design variables $\rho_i$ is equal to the sum of the implicit and explicit parts as follows:
\begin{align}
	\dv{f}{\rho_i} = \pdv{f}{\rho_i} + %
	\begin{Bmatrix}[1.5]
		\pdv{f}{\mathbf{u}} \\ \pdv{f}{\mathbf{x}} \\ \pdv{f}{\mathbf{w}}
	\end{Bmatrix}^T %
	\begin{Bmatrix}[1.5]
		\pdv{\mathbf{u}}{\rho_i} \\ \pdv{\mathbf{x}}{\rho_i} \\ \pdv{\mathbf{w}}{\rho_i}
	\end{Bmatrix}.
\end{align} 

The first term on the RHS can be calculated analytically. In order to calculate the derivatives of the state variables ($\mathbf{u}, \, \mathbf{x}, \, \textrm{and } \mathbf{w}$) w.r.t. the design variables $\rho_i$, the derivatives of the governing equations with respect to the design variables have to be used as in:
\begingroup
\allowdisplaybreaks
\begin{align}
	& \begin{Bmatrix}[1.5]
		\pdv{\mathbb{S}}{\rho_i} \\ \pdv{\mathbb{D}}{\rho_i} \\ \pdv{\mathbb{F}}{\rho_i}
	\end{Bmatrix} + %
	\begin{bmatrix}[1.5]
		\pdv{\mathbb{S}}{\mathbf{u}} & \pdv{\mathbb{S}}{\mathbf{x}} & \pdv{\mathbb{S}}{\mathbf{w}} \\
		\pdv{\mathbb{D}}{\mathbf{u}} & \pdv{\mathbb{D}}{\mathbf{x}} & \pdv{\mathbb{D}}{\mathbf{w}} \\
		\pdv{\mathbb{F}}{\mathbf{u}} & \pdv{\mathbb{F}}{\mathbf{x}} & \pdv{\mathbb{F}}{\mathbf{w}}
	\end{bmatrix} %
	\begin{Bmatrix}[1.5]
		\pdv{\mathbf{u}}{\rho_i} \\ \pdv{\mathbf{x}}{\rho_i} \\ \pdv{\mathbf{w}}{\rho_i}
	\end{Bmatrix} = \mathbf{0}, \\
	& \begin{Bmatrix}[1.5]
		\pdv{\mathbf{u}}{\rho_i} \\ \pdv{\mathbf{x}}{\rho_i} \\ \pdv{\mathbf{w}}{\rho_i}
	\end{Bmatrix} = %
	\begin{bmatrix}[1.5]
		\pdv{\mathbb{S}}{\mathbf{u}} & \pdv{\mathbb{S}}{\mathbf{x}} & \pdv{\mathbb{S}}{\mathbf{w}} \\
		\pdv{\mathbb{D}}{\mathbf{u}} & \pdv{\mathbb{D}}{\mathbf{x}} & \pdv{\mathbb{D}}{\mathbf{w}} \\
		\pdv{\mathbb{F}}{\mathbf{u}} & \pdv{\mathbb{F}}{\mathbf{x}} & \pdv{\mathbb{F}}{\mathbf{w}}
	\end{bmatrix}^{-1}
	\begin{Bmatrix}[1.5]
		\pdv{\mathbb{S}}{\rho_i} \\ \pdv{\mathbb{D}}{\rho_i} \\ \pdv{\mathbb{F}}{\rho_i}
	\end{Bmatrix}.
\end{align}
\endgroup

As for the details of the hyperelastic material model, we follow the excellent treatment in \citet[p.~492]{Bathe2014finite}. In fact, \(\mathbb{D}\) represents the residual vector at equilibrium as in:
\begin{align}
	& \mathbb{D} = {}^{ex}\mathbf{F} - {}^{in}\mathbf{F} = {}^t\mathbf{K} \; \Delta \mathbf{x} \simeq \mathbf{0}.
\end{align}

\noindent where ${}^{ex}\mathbf{F}$ is the external nodal force\textcolor{black}{s} vector, ${}^{in}\mathbf{F}$ is the internal nodal force\textcolor{black}{s} vector, ${}^t\mathbf{K}$ is the tangent stiffness matrix, and $\mathbf{x}$ is the incremental nodal displacement\textcolor{black}{s} vector (which is almost zero at equilibrium). This simply means that, \textit{at equilibrium}, the internal nodal forces due to material stresses are \textit{almost equal} to the external nodal forces within a certain tolerance. In fluid mesh deformation (be it a hyperelastic model or otherwise), there are prescribed displacements at the fluid-structure interface due to the structural deformations (i.e., non-homogeneous boundary conditions) and no prescribed\footnote{The term \textit{prescribed} is used here to describe both zero as well as non-zero prescribed displacements. Typically, zero displacements occur at the boundaries of the fluid mesh while non-zero ones occur at the fluid-structure interface due to the structural deformations.} external nodal forces. Luckily, the non-hom- ogeneous boundary conditions can be treated internally within ${}^{in}\mathbf{F}$ and ${}^t\mathbf{K}$ without the need to explicitly calculate the unknown external reactions corresponding to the prescribed displacements. This treatment of non-homogeneous boundary conditions is of interest in this context since it applies at the final equilibrium state as well. We detail the part pertaining to the tangent stiffness matrix as follows:
\begin{align}
	& {}^t\mathbf{K}_{\Omega_p, \Omega_{a \backslash p}} = \mathbf{0}, \label{eq:kt_treat1} \\
	& {}^t\mathbf{K}_{\Omega_p, \Omega_p} = \mathbf{I}. \label{eq:kt_treat2}
\end{align}

\noindent
where the subscript $\Omega_p$ designates the fluid mesh nodes with prescribed displacements (i.e., on the fluid-structure interface) and the subscript $\Omega_a$ designates all the fluid mesh nodes. As for the part that pertains to the internal nodal forces, it can be detailed as follows (though it isn't relevant to our discussion):
\begin{align}
	{{}^{in}\mathbf{F}}_{\Omega_p} & = L_f \, \mathbf{u}_{\Omega_p}  \\
	& (\textrm{at the start of each load incremental increase}), \nonumber \\
	{{}^{in}\mathbf{F}}_{\Omega_p} & = \mathbf{0} \\
	& (\textrm{at all the other Newton-Raphson iterations}). \nonumber
\end{align}

\noindent
where $L_f$ is a load factor that should sum up to 1 from all the load incremental steps used. A nice property of interest here is the relation between the internal nodal forces ${}^{in}\mathbf{F}$ and the tangent stiffness matrix ${}^t\mathbf{K}$ (see \citet[p.~493]{Bathe2014finite}):
\begin{align}
	\pdv{{}^{in}\mathbf{F}}{\mathbf{x}} = {}^t\mathbf{K}.
\end{align}

\noindent which will be utilized at the equilibrium state to calculate the required derivatives. For the derivative of \(\mathbb{D}\) with respect to the structural displacements $\mathbf{u}$:
\begin{align}
	\pdv{\mathbb{D}}{\mathbf{u}} = - \pdv{{}^{in}\mathbf{F}}{\mathbf{u}} = - \pdv{{}^{in}\mathbf{F}}{\mathbf{x}} \pdv{\mathbf{x}}{\mathbf{u}} = - {}^t\mathbf{K} \pdv{\mathbf{x}}{\mathbf{u}}.
\end{align}

\noindent
where $\mathbf{x}_{,\mathbf{u}}$ is a mapping that can be calculated from the relation linking the structural displacements to fluid mesh displacements. In this context, it makes sense to simply use conformal mesh at the fluid-structure interface so the same mesh nodes are shared between the solid and the fluid mesh. In this case, $\mathbf{x}_{,\mathbf{u}}$ would contain 1's at the shared nodes and 0's everywhere else and would be very cheap to compute. This relation can be stated as follows:
\begin{align}
	& \mathbf{x} = \mathbf{N} \mathbf{u}, \\
	& \mathbf{N}_{\Omega_p, \Omega_p} = \mathbf{I}, \\
	& \mathbf{N}_{\Omega_{a \backslash p}, \Omega_{a \backslash p}} = \mathbf{0}, \\
	& \pdv{\mathbf{x}}{\mathbf{u}}_{\Omega_p, \Omega_p} = \mathbf{I}, \\
	& \pdv{\mathbf{x}}{\mathbf{u}}_{\Omega_{a \backslash p}, \Omega_{a \backslash p}} = \mathbf{0}.
\end{align}

If non-conformal mesh is to be used, other techniques can be used to calculate this mapping (see \cite{Farhat1998b}). As for the derivative of \(\mathbb{D}\) with respect to the fluid mesh displacements $\mathbf{x}$:
\begin{align}
	& \pdv{\mathbb{D}}{\mathbf{x}} = \pdv{{}^{ex}\mathbf{F}}{\mathbf{x}} - \pdv{{}^{in}\mathbf{F}}{\mathbf{x}}.
\end{align}

Recall that there are no external nodal forces in fluid mesh deformation. Hence, the derivative can be simplified to: 
\begin{align}
	\pdv{\mathbb{D}}{\mathbf{x}} = - {}^t\mathbf{K}.
\end{align}

Finally, recall that the tangent stiffness matrix ${}^t\mathbf{K}$ is calculated using relations that can be found in any textbook on nonlinear finite element analysis, see \citet[p.~542]{Bathe2014finite}, then the non-homogeneous boundary conditions treatment in Eqs. \ref{eq:kt_treat1} and \ref{eq:kt_treat2} must be applied before proceeding with the governing equations' derivatives. Note that this sensitivity analysis is valid for any hyperelastic material model and not necessarily for the Yeoh model only. The only difference is in how the tangent matrix ${}^t\mathbf{K}$ is calculated.

\section{Conclusions}
\label{sec:conclusions}
In this work, we focused on mesh deformation techniques of structured, quadrilateral elements \textcolor{black}{to be considered} for implementation in topology optimization of fluid-structure interactions. We emphasized mesh quality rather than mesh admissibility as maintaining high mesh quality would inherently ensure its admissibility. As gradient-based optimization is the de facto approach in topology optimization, having smooth and continuous derivatives is a necessity in any mesh deformation technique. In our discussion of mesh quality, we paid more attention to shape changes (i.e., minimum skewness) rather than volume changes (i.e., minimum and maximum area changes). This is mainly because inverted elements and overlapping - which result in mesh inadmissibility - are mainly prevented by controlling shape rather than volume changes. In addition, the risk of element collapsing in structured, quadrilateral meshes is not as severe as in unstructured, triangular meshes.

We discussed and suggested some improvements to two legacy mesh deformation techniques; namely the spring analogy model and the linear elasticity model. In addition, we proposed a new mesh deformation technique based on the Yeoh hyperelasticity model. In terms of mesh quality, the hyperelastic model came first, then came the spring analogy model and the linear elasticity model respectively. In terms of the computational costs, the exact reverse order applies.

In all mesh deformation models discussed in this work, layered selective stiffening was applied to consecutive layers of elements adjacent to the fluid-structure interface where the bulk of mesh distortion is located. Stiffening was achieved through multiplying the material constant (the elastic modulus \textcolor{black}{in the legacy models} or the material parameter $A_{20}$ \textcolor{black}{in the hyperelastic model}) by a factor. The benefit of layered selective stiffening is twofold. Increasing the minimum skewness of the whole mesh and moving the distorted elements away from critical fluid flow regions. \textcolor{black}{Two comments about layered selective stiffening are in order:
\begin{enumerate}
	\item It's true that in this work we assume the fluid-structure interface to be known apriori, hence the identification of layers of finite elements to be stiffened is straightforward. However, in a true TOFSI with large deformations the fluid-structure interface would be continuously evolving. In order to apply this layered selective stiffening to an implicit fluid-structure interface, we propose the following. In most cases, we conjecture that the majority of the mesh deformation would be attained in the first solution of the governing equations. More so if the initial structural density is the volume fraction applied uniformly to the whole design space which is the typical approach in density-based TO. The stiffening factors could then be applied as a function of the elemental density (i.e., design variables) so gray elements would have some degree of stiffening while they transition to either complete solid or fluid elements. In most density-based TO problems, a rough topology is usually reached early in the optimization process while the remaining iterations are used to stabilize and tune this rough result. Hence, layered selective stiffening would be applied more accurately once this rough design is reached.
	\item The data presented in this work could potentially be helpful in providing an educated guess on the value of the stiffening factors to be used with each mesh deformation model. Since the optimal stiffening factors seem to vary widely even for the same initial geometry under different deformation modes, it would take some trials and errors to approach these optimal factors for different problems. Nonetheless, as long as the stiffening factors are decreased as layers move away from the fluid-structure interface, the mesh quality would still improve even if not optimally.
\end{enumerate}
A final noteworthy remark is that the computational cost of deforming the fluid mesh could be greatly reduced by localizing the mesh deformation to the area surrounding the moving structure where the bulk of the mesh deformation is located.}

\vspace{1em}
\noindent \textbf{\large Declarations}

\noindent \textbf{Conflict of interest} The authors declare they have no conflict of interest.

\noindent \textbf{Replication of results} The authors have included all dimensions and material parameters needed to replicate the results.

\bibliographystyle{spbasic}
\pdfbookmark[0]{References}{References}
\bibliography{library}

\end{document}